\documentclass[11pt]{amsart}%
\usepackage{xypic}
\xyoption{all}
\usepackage{amsthm,euscript} 
\usepackage{paralist} \usepackage{amsopn}
\usepackage{amsfonts}
\usepackage{amsmath}
\usepackage{latexsym}
\usepackage{amscd}
\usepackage{amssymb}
\usepackage{amsmath}
\usepackage{fancyhdr}
\usepackage{setspace}

\usepackage[colorlinks=true, pdfstartview=FitV, linkcolor=blue,
citecolor=blue,
urlcolor=blue,breaklinks=true]{hyperref}

\long\def\ignore#1\endignore{\relax}

\newcommand{\sm}{\smallskip}
\newcommand{\ms}{\medskip}

\usepackage{fancyhdr}

\swapnumbers
\newtheorem{theorem}[subsection]{Theorem}

\newtheorem{proposition}[subsection]{Proposition}
\newtheorem{prop}[subsection]{Proposition}
\newtheorem{lemma}[subsection]{Lemma}
\newtheorem{Lemma}[subsection]{Lemma}
\newtheorem{corollary}[subsection]{Corollary}

\theoremstyle{definition}
\newtheorem{defn}[subsection]{Definition}
\newtheorem{remark}[subsection]{Remark}

\newtheorem{example}[subsection]{Example}

\numberwithin{equation}{section} \allowdisplaybreaks

\numberwithin{equation}{subsection}

\def\gg{\mathfrak g}

\def\SS{\mathfrak S}

\def\mmm{\mathfrak m}
\def\nn{\mathfrak n}

\def\pp{\mathfrak p}
\def\Ralg{R\text{--}\mathbf{alg}} 
\def\et{\acute et}

\newcommand{\bbP}{{\mathbb P}}
\newcommand{\bbQ}{{\mathbb Q}}
\newcommand\NN{\mathbb{N}}

\newcommand{\Aut}{\operatorname{Aut}}
 \DeclareMathOperator{\ad}{ad}

\newcommand{\End}{\operatorname{End}}

\newcommand{\Ctd}{\operatorname{Ctd}}

\newcommand{\Lie}{\operatorname{Lie}}

\newcommand{\GL}{{\operatorname{GL}}}

\newcommand{\Pic}{{\operatorname{Pic}}}
\newcommand{\simlgr}{\buildrel \sim \over \lra}
\newcommand{\PGL}{{\operatorname{PGL}}}

\newcommand{\Int}{{\operatorname{Int}}}

 \newcommand{\bAut}{\rm \bf Aut}

\newcommand{\diag}{\operatorname{diag}}

\newcommand{\lra}{\longrightarrow}

\newcommand{\Mat}{{\operatorname{M}}}

\newcommand{\tr}{\operatorname{tr}}

\newcommand{\cal}{\mathcal}

\def\bF{\text{\rm \bf F}}
\def\bH{\text{\rm \bf H}}
\def\bG{\text{\rm \bf G}}

\def\bAut{\mathbf{Aut}} 

\def\al{\alpha}

\def\si{\sigma}

\newcommand\co{\colon} 
\newcommand\wti{\widetilde}
\newcommand\ideal{\triangleleft}
\newcommand\ot{\otimes}

\newcommand\op{^{\rm op}}

\DeclareMathOperator{\Span}{span}

\DeclareMathOperator{\supp}{supp}

\newcommand\de{\delta}
\newcommand\eps{\varepsilon}
\newcommand\ga{\gamma}

\newcommand\io{\iota}
\newcommand\ka{\kappa}
\newcommand\la{\lambda}
\newcommand\La{\Lambda} 

\newcommand\ta{\tau}
\newcommand\vphi{\varphi} 
\newcommand\ze{\zeta}
\newcommand\om{\omega}

\newcommand\boldeps{{\boldsymbol{\varepsilon}}}
\newcommand\beps{\boldeps}

\newcommand{\euZ}{\EuScript{Z}}

\newcommand\scA{\mathcal{A}}
\newcommand\scB{\mathcal{B}}
\newcommand\scC{\mathcal{C}}
\newcommand\scD{\mathcal{D}}

\newcommand\scF{\mathcal{F}}

\newcommand\scL{\mathcal{L}}

\newcommand\scO{\mathcal{O}}

\newcommand\scZ{\mathcal{Z}}

\newcommand\frh{\ensuremath{\mathfrak{h}}} \newcommand\h{\frh}
\newcommand\lsl{\ensuremath{\mathfrak{sl}}}
\newcommand\gl{\ensuremath{\mathfrak{gl}}}

\newcommand\ZZ{\mathbb{Z}}

\newcommand\sfB{{\sf B}}

\newcommand{\bGL}{\mathbf{GL}}

\newcommand{\bPGL}{\mathbf {PGL}}

\begin{document}

\title[]{On conjugacy of Cartan subalgebras
in non-fgc Lie tori}

\author{V. Chernousov}
\address{Department of Mathematics, University of Alberta,
    Edmonton, Alberta T6G 2G1, Canada}
\thanks{ V. Chernousov was partially supported by the Canada Research
Chairs Program and an NSERC research grant} \email{chernous@math.ualberta.ca}

\author{E. Neher}
\address{Department of Mathematics and Statistics, University of Ottawa,
    Ottawa, Ontario K1N 6N5, Canada}
\thanks{E.~Neher was partially supported by a Discovery grant from NSERC}
\email{neher@uottawa.ca}

\author{A. Pianzola}
\address{Department of Mathematics, University of Alberta,
    Edmonton, Alberta T6G 2G1, Canada.
    \newline
 \indent Centro de Altos Estudios en Ciencia Exactas, Avenida de Mayo 866, (1084) Buenos Aires, Argentina.}

\thanks{A. Pianzola wishes to thank NSERC and CONICET for their
continuous support}\email{a.pianzola@gmail.ca}

\thanks{The first two authors would like to thank CAECE for their hospitality while most of this paper was written.}

\begin{abstract} We establish the conjugacy of Cartan subalgebras for generic Lie tori ``of type A". This is the only conjugacy problem  of Lie tori related to Extended Affine Lie Algebras that remained open.

\end{abstract}

\maketitle

\hfill
{\em To H. Abels on his 75th birthday}

\section*{Introduction}

Extended Affine Lie Algebras (EALAs for short) are a rich class of Lie
algebras  that were first conceived by the physicists R.~H{\o}egh-Krohn and
B.~Torresani and then brought to the attention of mathematicians by P.
Slodowy (the reader should look  at \cite{n:persp} and \cite{n:eala-summ} for
a comprehensive review of basic EALA  theory and references. The original
mathematical formulation is to be found in \cite{AABGP}). An EALA has an
invariant non-negative integer attached called its {\it nullity}. In nullity
0 EALAs are nothing but the finite dimensional simple Lie algebras, while in
nullity 1 EALAs are the celebrated affine Kac-Moody Lie algebras (we assume
for simplicity in this Introduction that our base field is the complex
numbers). Roughly speaking an EALA $E$ is constructed from a class of Lie
algebras called Lie tori by taking a central extension and adding a suitable
space of derivations. In the case of the affine algebras, for example, the
Lie torus is a loop algebra $L$ based on a finite dimensional simple Lie
algebra $\gg$. Note that $L$ is naturally a Lie algebra over the Laurent
polynomial ring $\mathbb{C}[t^{\pm 1}]$). This ring is the centroid of
$L.$The central extension of $L$ is the universal one (which is
one-dimensional). The space of derivations is also one-dimensional and
corresponds to the degree derivation $t(d/dt)$.

An EALA, by definition, comes equipped with a so-called  Cartan subalgebra
(just like the affine algebras do, but unlike the finite dimensional simple
Lie algebras; see \S\,\ref{lsl_2} for details). In the setting of EALAs, a Cartan 
subalgebra is the same as a self-centralizing ad-diagonalizable subalgebra, 
as defined in \S\,\ref{def:mad}.
With respect to the given Cartan subalgebra the EALA admits a
root space decomposition. The structure of the resulting ``root system" plays
a fundamental role in understanding the structure as well as the
representation theory of the given EALA. It is obvious that all of this would
be of little use (or mathematically unnatural) if the nature of the root
system was to depend on the choice of Cartan subalgebra. The most elegant way
of dealing with this problem is by establishing ``Conjugacy", i.e.,  by
showing that all Cartan subalgebras are conjugate under the action of the
group of automorphisms of the EALA (in all cases it is sufficient to use a
precise subgroup of the full group of automorphisms. Conjugacy in the finite
dimensional case,  in the spirit of the present work,  is due to Chevalley.
For the affine algebras the result is due to Peterson and Kac \cite{PK}). For
almost all EALAs (see below) conjugacy, hence the invariance of  the root
system, has been established  in \cite{CGP} (for Lie tori) and \cite{CNPY}
(for the full EALAs). One case, the so called non-fgc case (see \S2 for definitions), remained open.
The purpose of this paper is to establish conjugacy for the Lie tori (the
analogue of \cite{CGP} in the non-fgc case) underlaying this remaining family
of EALAs.

The centroid of a Lie torus $L$ is always a Laurent polynomial ring $R$ in finitely many variables. In all cases but one, $L$ is an $R$-module of finite type. This is the fgc case (where fgc stands for finitely generated over the centroid). When this does not happen, the non-fgc case, the nature of $L$ is perfectly understood:
$$ L = \lsl_\ell(Q)$$
where $Q$ is a quantum torus with at least one generic entry. We remind the reader that by definition $Q$ is the complex  unital associative algebra presented by generators $x_1, \ldots, x_n , x_1^{-1}, \ldots,
x_n^{-1}$ and relations \[ x_i x_i^{-1} = 1_Q = x_i^{-1} x_i, \quad x_i x_j =
q_{ij} x_j x_i. \] where the $q_{ij}$ are non-zero complex numbers.  That ``one entry be generic" means that one of the $q_{ij}$ cannot be a root of unity. The centre $\euZ(Q)$ of $Q$ is always a Laurent polynomial ring. The fgc condition on the Lie algebra $\lsl_\ell(Q)$ is equivalent to $Q$ being a $\euZ(Q)$-module of finite type.

We can now state our
\bigskip

{\bf Main Theorem.} {\it Let $ f \co \lsl_{\ell'}(Q') \to \lsl_\ell(Q)$ be an
isomorphism of non-fgc Lie tori. Then $\ell = \ell'$ and if\/ $\h'$ and $\h$
denote the given Cartan subalgebras of\/ $ \lsl_{\ell'}(Q')$ and\/
$\lsl_\ell(Q)$ respectively,\footnote{As we have mentioned already, a
distinguished ``Cartan subalgebra" is part of the definition of a Lie torus.}
then $f(\h')$ and $\h$ are conjugate under an automorphism of
$\lsl_\ell(Q).$}

\bigskip
In the fgc case the proof of conjugacy is  naturally divided into two steps.
One first establishes conjucacy at the level of Lie tori, and then extends
this conjugacy ``downstairs" to the full EALA. The fgc condition allows the
Lie tori to be viewed as simple Lie algebras (in the sense of \cite{SGA3})
over Laurent polynomial rings. Conjugacy in the fgc case makes heavy use of
the powerful methods of \cite{SGA3} and Bruhat-Tits theory. None of this is
possible in the non-fgc case. New methods/ideas are needed. The crucial
ingredient that we develop to deal with this new situation is a method that
we call ``specialization". The idea, roughly speaking, is to create a subring
$R$ of $\mathbb{C}$ with the property that

(i) our non-fgc Lie torus ``exists" over $R,$

(ii)there exists a maximal ideal $\mmm$ of $R$ which after base change
(reduction modulo $\mmm$) yields an fgc Lie torus.

The catch is that the field $R/\mmm$ is of positive characteristic! One does
not even have a suitable definition of Lie tori in positive characteristic.
Yet the resulting object and its group of automorphisms is explicit enough
that we can establish conjugacy for them. The specialization method is
invoked once again to show that conjugacy holds before the reduction modulo
$m$.

\bigskip
{\em Notation.} Throughout, $R$ is a commutative unital ring which often occurs as the
base ring of some algebraic structure being considered; $F$ is an arbitrary field, and $k$ often
denotes a field of characteristic $0$. An $R$--algebra is an arbitrary algebra over $R$ (in particular not
necessarily associative or a Lie algebra). Group schemes are usually denoted
with bold letters. For example, $\bPGL_{m,R}$ denoted the $R$-group scheme of automorphisms of the
(associative and unital) $R$-algebra $M_m(R).$

\bigskip

 {\em Structure of the paper.} In Section 1 we collect some basic results about centroids of (arbitrary) algebras. Particular attention is devoted  to the case of quantum tori. Section 2 looks at the structure of the Lie algebra $\lsl_\ell(Q)$ where $Q$ is a quantum torus. The definition and basic properties of maximal abelian diagonalizable (MAD) subalgebra are also given in this section (these are the subalgebras that play the role of the Cartan subalgebras in EALA theory). Section 3 is devoted to a detailed study of the group scheme $\bPGL_{\mathcal{A}}$ where $\mathcal{A}$ is an Azumaya algebra over a ring $R$, and the connection between this $R$-group scheme and  the $R$-group scheme of automorphisms of $\lsl_\ell(\mathcal{A}).$ Section 4 presents a detailed analysis of the automorphism group of the Lie algebra $\lsl_\ell(Q)$ when $Q$ is an fgc quantum torus. Section 5 develops the method that we called ``specialization" mentioned above. This is the key that allows us to deal with the non-fgc case by translating the problem into an fgc question, but now over fields of positive characteristic. Section 6 presents a collection of preliminary results to be used in the proof of the main Theorem, which is given in Section 7.
\section{Some results on centroids and quantum tori}\label{sec:revQ}

\subsection{Centroids and base change}\label{cbc}
 Let $R$ be a commutative ring and let $\scA$
be an arbitrary $R$--algebra\footnote{It will not be sufficient in the
following to consider only algebras over fields}. Recall that the derived subalgebra of $\scA$ is the additive subgroup of $\scA$ generated by all products $ab$ with $a,b \in \scA$. It is trivial to see that this group is indeed an $R$-subalgebra of $\scA$. The algebra $\scA$ is called {\it perfect} if it equals to its derived algebra.  Note that any unital algebra is perfect.

 A crucial object for our work is the {\em centroid\/} $\Ctd(\scA)$ of the algebra $\scA$. Recall that
\[ \Ctd(\scA) = \{ \chi \in\End_R(\scA) : \chi(a_1 a_2) =  \chi(a_1) a_2 =
    a_1 \chi(a_2) \text{ for all } a_i \in \scA \}.
\]
Clearly $\Ctd(\scA)$ a unital commutative (if $ \scA$ is perfect) subalgebra of the associative $R$--algebra $\End_R
(\scA)$. It is obvious that we can consider $\scA$ as an algebra over $\scC =
\Ctd(\scA)$ -- it will be denoted $\scA_{(\scC)}$.\footnote{We use
$\scA_{(\scC)}$ instead of $\scA_\scC$ since the latter usually denotes base
change.} We will say that $\scA$ is {\em fgc\/} if $\scA_{(\scC)}$ is a
finitely generated $\scC$--module.

More generally, if $S\in \Ralg$, i.e., $S$ is a unital associative
commutative $R$--algebra, and if $\rho \co S \to \Ctd(\scA)$ is a unital
algebra homomorphism, $\scA$ becomes an $S$--algebra by defining $ s \cdot a
= \rho(s)(a)$ for $s\in S$ and $a\in \scA$. We will denote the algebra
obtained in this way by $\scA_{(\rho)}$ or $\scA_{(S)}$ if $\rho$ is clear
from the context. \sm

{\bf Example 1.} Assume that $f\co \scA' \simlgr \scA$ is an isomorphism of
perfect $R$--algebras. It is then easily seen (and well known) that
\begin{equation} \label{cbc0}
   \Ctd(f) \co \Ctd(\scA) \simlgr \Ctd(\scA'), \quad \chi \mapsto
         f ^{-1} \circ \chi \circ f
\end{equation}
is an isomorphism of $R$--algebras. Since $\scC = \Ctd(\scA) \in \Ralg$, we
can use $\Ctd(f)$ to make $\scA'$ a $\scC$--algebra: $\chi \cdot a' =
\big(\Ctd(f) (\chi)\big)(a') = ( f^{-1} \circ \chi \circ f)(a')$. Then
\begin{equation} \label{cbc1}
f \co \scA'_{(\scC)} \to \scA_{(\scC)} \; \hbox{\em   is
$\scC$--linear.}
\end{equation}
Indeed,
\[f(\chi \cdot a') = f \big( (f^{-1} \circ \chi \circ f) (a')\big) =
\chi\big( f(a')\big) = \chi \cdot f(a')\] for $\chi \in \scC$ and $a'\in
\scA'$. \sm

{\bf Example 2 (centre).} Let $\scA$ be a unital associative $R$--algebra. As
usual, $[a,b] = ab-ba$ for $a,b,\in \scA$ denotes the (Lie) commutator. The
{\em centre\/} $\euZ(\scA)= \euZ $ consists of all $z\in \scA$ satisfying
$[z,a] = 0$. One easily checks that the left multiplication by a central
element is an isomorphism of $R$--algebras:
\begin{equation} \label{cbc2}
     \euZ(\scA) \simlgr \Ctd(\scA), \quad z \mapsto L_z.
\end{equation} Hence we can (and will)
consider $\scA$ as a $\euZ$--algebra, denoted $\scA_{(\euZ)}$. \sm

{\bf Example 3.} Let $\scA$ be as above and let $\scL = \lsl_\ell(\scA)$ be
the special linear Lie $R$--algebra introduced in \ref{def:lsl}. If  $\euZ = \euZ(\scA)$, we have an
obvious $R$--algebra homomorphism
\begin{equation*}
  \label{cbc22} \ze \co \euZ \to \Ctd(\scL), \quad \ze(z) = \big(
      (x_{ij}) \mapsto (zx_{ij}) \big)
\end{equation*}
for $z\in \euZ$ and $x=(x_{ij}) \in \scL$. We will show in Lemma~\ref{ctd-sl}
that $\ze$ is an isomorphism if $\ell \ge 2$ and $\frac{1}{2}\in R$. It is
easily seen that
\begin{equation}
  \label{cbc3}
     \lsl_\ell(\scA)_{(\scZ)} = \lsl_\ell(\scA_{(\scZ)} ) .
\end{equation}

\subsection{Some properties of quantum tori}\label{qua-pro}
We list some properties of quantum tori that we will use. Throughout, $F$
is a field of arbitrary characteristic, and $\La$ is a free abelian group of
rank $n$. \sm

\begin{inparaenum}[(a)]
\item \label{qua-pro1}(Definitions) By definition, a {\em quantum torus (with
    grading  group $\La$)\/} is an associative  unital $\La$-graded
    $F$-algebra   $Q= \bigoplus_{\la \in \La} Q^\la$ such
    that $\dim Q^\la = 1$ for all $\la \in \La$ and that every $0 \ne a \in Q^\la$ is
    invertible.

After fixing a basis $\beps= (\eps_i)$ of $\La$, we can choose $0 \ne x_i \in
Q^{\eps_i}$ and then get a quantum matrix $q=(q_{ij}) \in \Mat_n(F)$ defined
by $x_i x_j = q_{ij} x_j x_i$. Then, using $x_i^{-1} = $ the inverse of
$x_i$,  we define $x^\la = x_1^{\ell_1} \cdots x_n^{\ell_n}$ for $\la =
\ell_1 \eps_1 + \cdots + \ell_n \eps_n \in \La$:
\begin{equation}
  \label{qua-pro11} Q = \textstyle \bigoplus_{\la \in \La} F x^\la.
\end{equation}
One can then also realize a quantum torus as the unital associative
$F$-algebra presented by generators $x_1, \ldots, x_n , x_1^{-1}, \ldots,
x_n^{-1}$ and relations \[ x_i x_i^{-1} = 1_Q = x_i^{-1} x_i, \quad x_i x_j =
q_{ij} x_j x_i. \] We will refer to this view of $Q$ as a {\em
coordinatization}.

We point out that such a presentation is not unique: it depends on the chosen
$\ZZ$--basis $\beps$ of $\La$. In other words, for any integral matrix
$A=(a_{ij})\in \GL_n(\mathbb{Z})$ the set
$\tilde{x}=\{\tilde{x}_1,\ldots,\tilde{x}_n\}$ of invertible elements in $Q$,
defined by
\[
   \tilde{x}_1=x_1^{a_{11}}\cdots x_n^{a_{1n}},
\ldots,\tilde{x}_n=x_1^{a_{n1}}\cdots x_n^{a_{nn}},
 \]
also generates $Q$ and the associated quantum matrix
$\tilde{q}=(\tilde{g}_{ij})$ is given by $ \tilde{q}_{ij}=\prod_{s,t}
q_{st}^{a_{is}a_{jt}} $. \sm

\item\label{qua-pro2} The centre of $Q$ is a $\La$-graded
    subalgebra,
\[
 \euZ(Q) = \textstyle \bigoplus_{\xi  \in \Xi} Q^\xi
\]
where $\Xi$ is the so-called {\em central grading group\/}:
\[\Xi = \{ \la \in
\La : Q^\la \subset \euZ(Q)\}.\] This is a free abelian group of rank $z\le
n$. Hence $\euZ(Q)$ is a Laurent polynomial ring in $z$ variables, which we
may take as $t_1, \ldots, t_z$ (these can be taken to be of the form $x^\la$ for suitable
$\la$s). \sm

\item\label{qua-pro3}  We define
\[[Q,Q] = \Span_F \{[a,b]: a,b\in Q\},\] a
    graded subspace of $Q$. One knows (see e.g.\ \cite[Prop.~2.44(iii)]{bgk}
    for $F=\mathbb{C}$ or \cite[(3.3.2)]{ny} in general)
\begin{equation}\label{qua-pro33}
   Q = \euZ(Q) \oplus [Q,Q]. \end{equation}

\item\label{qua-pro4} $Q$ is a domain: $ab=0$ for $a,b\in Q$ implies that
    $a=0$ or $b=0$, whence a nondegenerate and thus prime associative
    $F$-algebra. This implies that $Q$ is connected: the only idempotents in $Q$
    are $0$ and $1_Q$.
\sm

\item\label{qua-pro-unit} An element $u$ of $Q$ is invertible if and only if
    $0\ne u \in Q^\la$ for some $\la \in \La$. \sm

\item \label{qua-pro-e} The grading properties of a quantum torus $Q$ show
    that     $Q$    is    fgc    in the sense of \ref{cbc} if and only if $\Xi$  has  finite index in $\La$.
    Equivalently, for some (hence all) coordinatization all entries $q_{ij}$ of the quantum
    matrix $q$ have finite order. If this holds, then for every coordinatization the
    $q_{ij}$ have finite order.

\item \label{qua-pro5} We let $\euZ = \euZ(Q)$, and denote by $\wti \euZ$ the
    quotient field
    of $\euZ$, a rational function field. The $\wti \euZ$-algebra \[ \wti Q  = Q
    \ot_\euZ \wti \euZ\] is called the {\em central closure of $Q$}. It has
    the following properties: $\wti Q$ is a central $\wti \euZ$-algebra,
    $\wti Q$ is a domain (since $Q$ is a domain), and $Q$ embeds into $\wti Q$. In particular,
    if $Q$ is fgc, then $\wti Q$ is a finite-dimensional central domain over $\wti \euZ$,
    whence a central division $\wti \euZ$-algebra.

\item \label{qua-protr} ({\em Trace, $\ZZ$-grading, degree}) As in
    \eqref{qua-pro1} we fix  a   basis $\boldeps =(\eps_1, \ldots, \eps_n)$ of $\La$ and define the
    {\em $\boldeps$-trace\/} as $\tr_{\boldeps}(\la) = \sum_i \la_i$ for $\la =
    \sum_i \la_i \eps_i$ with $\la_i \in \ZZ$. Since
$\tr_\beps \co \La \to \ZZ$
is a group homomorphism, the $\La$-grading of $Q$ can be made into a
$\ZZ$-grading
\begin{equation} \label{qua-protr0}
  Q =\textstyle  \bigoplus_{n \in \ZZ} Q_{(n,\beps)} \quad
       \text{with} \quad Q_{(n,\beps)} = \bigoplus_{\tr_\beps(\la) = n} Q^\la.
\end{equation}
Every $0\ne q\in Q$ can be uniquely written as $q = \sum_{n\le m} q_{(n)}$
with $q_{(n)} \in Q_{(n, \beps)}$ and $q_{(m)}\ne 0$. We call $m$ the {\em
$\beps$-degree of q\/} and denote it by $\deg_\beps q$. We put $\deg_\beps 0
= - \infty$.

In the following we may suppress the dependance on $\boldeps$ and just speak
of the trace. Analogously for the $\ZZ$-grading
\eqref{qua-protr0}.
\end{inparaenum}

\begin{Lemma}\label{deg-prop} Let $Q$ be a quantum torus over the field $F$, and let
$q,q_1, q_2 \in Q$.
\begin{inparaenum}[\rm (a)]

\item \label{deg-prop-a} For any $\ZZ$-basis $\beps$ of $\La$ we have
\begin{align}
     \label{deg-prop2}
  \deg_\beps (q_1 + q_2) &\le \max \{ \deg_\beps(q_1), \deg_\beps(q_2)\}, \\
 \label{deg-prop3}
      \deg_\beps (sq) &= \deg_\beps (q) \quad \text{for $0\ne s\in F$},
 \\ \label{deg-prop1}
 \deg_\beps(q_1 q_2) &= \deg_\beps(q_1) + \deg_\beps(q_2)
\end{align}
with the obvious rules in case one of $q,q_1$ or $q_2 = 0$. \sm

\item \label{deg-prop-b} For $0 \ne q \in Q$ we have $q\in Q^0 = F 1_Q \iff
    \deg_\beps(q) =     0$     for     all     $\ZZ$-bases $\beps$ of $\La$.
\end{inparaenum} \end{Lemma}

\begin{proof} \eqref{deg-prop-a} All formulas are easily verified in case one of $q$, $q_1$ or $q_2 = 0$ (using the convenience that $0$ is of degree $-\infty$).
Also, \eqref{deg-prop2} and \eqref{deg-prop3} follow immediately from the
definition. For \eqref{deg-prop1} and $q_1, q_2 \ne 0$ we have $q_1 q_2 \ne
0$ since $Q$ is a domain. We write $q_1 = \sum_{n \le m} q_{(n)}$ with
$\deg_\beps q_1 = m$ and $q_2 = \sum_{n \le p } q'_{(n)}$ with $\deg_\beps
(q_2) = p$. Since $Q = \bigoplus_{n\in \ZZ} Q_{(n,\beps)}$ is a
$\ZZ$-grading, it follows that $q_1 q_2 = q_{(m)} q'_{(p)} + $ terms of lower
degree in the $\ZZ$-grading with respect to $\beps$.

\eqref{deg-prop-b} Let $0\ne q \in Q$. If $q\in Q^0$ then $\deg_\beps(q) = 0$
because $\tr_\beps(q) = 0$ for any $\ZZ$-basis $\beps$ of $\La$. Conversely,
we can assume $q\ne 0$ and write $q = \sum_{\la \in \supp(q)} q^\la$ with
$q^\la \in Q^\la$ and $\supp(q) = \{ \la \in \La: q^\la \ne 0 \}$. Then
$\tr_\beps(\la) \le 0$ for every $\la \in \supp(q)$ by assumption. But since
both $\beps$ and $-\beps$ are $\ZZ$-bases of $\La$ and $\tr_{-\beps} = -
\tr_\beps$ we get $\tr_\beps(\la) = 0$ for every $\la \in \supp(q)$. For a
fixed $\beps$ and $\la = \sum_{i=1}^n \ell_i \eps_i \in \supp(q)$ we
therefore have $\sum_{i=1}^n \ell_i = 0$. Our claim obviously holds if $\La
\cong \ZZ$. Thus we can assume that $\La$ has rank at least $2$. With respect
to the $\ZZ$-basis $\beps' = (\eps_1 + 2 \eps_2, \eps_2, \eps_3, \ldots)$ we
have $\la = \ell_1(\eps_1 + 2 \eps_2) + (\ell_2 - 2\ell_1) \eps_2 + \cdots$.
Thus $0 = \sum_i \ell_i = \ell_1 + (\ell_2 - 2\ell_1) + \sum_{i\ge 3} \ell_i
= 0$, and $\ell_1=0$ follows. Similarly, all $\ell_i= 0$, i.e., $\la = 0$.
\end{proof}

In the following lemma we will describe certain $F$-diagonalizable
endomorphisms $\phi $ of a quantum torus $Q$ over $F$. The term {\em
$F$-diagonalizable\/} means of course that there exists an $F$-basis of the
$F$-vector space $Q$ consisting of eigenvectors of $\phi$ with eigenvalues in
$F$.

\begin{Lemma}\label{diag-q}
Let $Q = \bigoplus_{\la \in \La} Q^\la$ be a quantum torus over the field
$F$, and let $d\in Q$. \begin{inparaenum}[\rm (a)]

\item \label{diag-q-a}  If $dq = \om q$ for some $0\ne q\in Q$ and
    $\om     \in   F$, then $d\in Q^0$. In particular, the left multiplication $L_d$ for $d\in Q$ is
    $F$-diagonalizable if and only if $d\in F 1_Q= Q^0$.

\item \label{diag-q-b}  If $[d,q] = \om q$ for some $0 \ne \om \in
    F$,      then    $q=0$.    In    particular, the endomorphism $\ad d\in \End_F(Q)$, defined by $(\ad d)(q)
    = [d,q]$, is $F$-diagonalizable if and only if $d\in \euZ(Q)$.
\end{inparaenum} \end{Lemma}

\begin{proof} (a) follows from the fact that $Q$ is a domain and $q\ne 0$.

(b) Suppose $0\ne q$. Then for any $\ZZ$-basis $\beps$ of $\La$ we get,
using the formulas of Lemma~\ref{deg-prop}\eqref{deg-prop-a},
\begin{align*}
  \deg_\beps(dq) &= \deg_\beps(d) + \deg_\beps(q)
                    = \deg_\beps(d) + \deg_\beps(-q) = \deg_\beps(-qd), \\
  \deg_\beps(q) &= \deg_\beps(\om q) = \deg_\beps(dq-qd)
    \le \max\{\deg_\beps(dq), \deg_\beps(-qd)\} \\ &  = \deg_\beps(dq)
           = \deg_\beps(d) + \deg_\beps(q)
\end{align*}
whence $\deg_\beps(d) \ge 0$. But since $-\beps$ is also a basis of $\La$, we
in fact have $\deg_\beps(d) = 0$. Thus, by
Lemma~\ref{deg-prop}\eqref{deg-prop-b}, we have $d\in Q^0$.
But then $[d,q]=0$ yields a contradiction. This shows $q=0$, and also that
$\ad d$ does not have a non-zero eigenvalue. In particular, if $\ad d$ is
$F$-diagonalizable, necessarily $\ad d=0$, i.e., $d\in \euZ(Q)$.
\end{proof}

In the remainder of this section we present some results which are special
for fgc quantum tori. \sm

\subsection{Canonical presentation}\label{canpre}
Let $Q$ be a quantum torus over the field $F$, coordinatized as $Q=
\bigoplus_{\la \in \La} F x^\la$ as in \eqref{qua-pro11}, and let
$q=(q_{ij})$ be the associated quantum matrix. We will say that the
coordinatization (or presentation) is {\it canonical} if all entries of the
quantum matrix $q$ outside of its diagonal blocks of size $2\times 2$ are
equal to $1$. Equivalently, for every $i\geq 1$ the generators
$x_{2i-1},x_{2i}$ of $Q$ commute with all other generators $x_j$ where
$j\not=2i-1,2i$.

The following lemma is stated without proof in \cite[Remark~7.2]{CP}. A proof
is given in \cite[Thm.~4.5]{N}.

\begin{lemma}\label{capr} Any fgc quantum torus $Q$ has a
canonical presentation.
\end{lemma}

\begin{example}[{\bf Quantum $2$-tori}]\label{n=2} 
Let $Q$ be a quantum torus
whose grading group has rank $2$. Hence  $Q$ is generated by two elements,
say $x_1$, $x_2$, and the corresponding quantum matrix has the form
\[
     q=\left(\begin{array}{cc}
1 & q_{12}\\
q_{12}^{-1} & 1
\end{array}
\right)
\]
where $q_{12} \in F^\times$. By definition, this presentation of $Q$ is
canonical. The algebra $Q$ is fgc if and only if $q_{12}$ is a root of unity,
say primitive of degree $\ell$. Let us assume this in the following.The $\euZ$--algebra $Q_{(\euZ)}$ will be
called a {\it symbol algebra} of degree $\ell$.
Its centre $\euZ$ is a Laurent
polynomial ring (\ref{qua-pro}\eqref{qua-pro2})
$$\euZ=k[t_1^{\pm 1},t_2^{\pm 1}]$$
where $t_1=x_1^\ell$, $t_2=x_2^\ell$ (observe that $\ell$ is independent of the
coordinatization, since $\La/\Xi = \ell^2$). We will usually denote this
$\euZ$--algebra by $(t_1,t_2)_{\euZ,\, q_{12}}$ or simply $(t_1,t_2)$ if
there is no risk of confusion. Note that $Q_{(\euZ)}$ has order $\ell$ in the
Brauer group ${\rm Br}(\euZ)$.

If $\ell \in F^\times$ the subalgebra $E= \euZ[x_1^{\pm 1}]$ is a maximal
(abelian)  \'etale subalgebra of $Q_{(\euZ)}$.
\end{example}

\begin{example}\label{tensorproduct}  Let $Q= \bigoplus_{\la \in \La} F x^\la$ be
an fgc quantum torus over $F$. By Lemma~\ref{capr} we can assume that $Q$ is
canonically presented, say with quantum matrix $q=(q_{ij})$. Up to
re-numbering (=re-coordinatization), we may assume that
$q_{12},q_{34},\ldots,q_{2s-1,2s}\not=1$, but $q_{2i+1,2i+2}=1$ for all
$i\geq s$. Then $Q$ admits a decomposition
\[
    Q_{(\euZ)}=Q_{1,\euZ}\otimes_\euZ \cdots\otimes_\euZ Q_{s,\euZ}
\]
where the $Q_{i,\euZ}=(t_{2i-1},t_{2i})$ are the symbol algebras in degree
$\ell_i$ corresponding to the nontrivial diagonal blocks of $q$ of size $2\times
2$, i.e. to the block diagonal sub-matrices
$$
\left(
\begin{array}{cc}
1 & q_{2i-1,2i}\\
q_{2i,2i-1} & 1
\end{array}
\right)
$$
of $q$ where $i\leq s$. Here $t_{2i-1}=x_{2i-1}^{\ell_i},t_{2i}=x_{2i}^{\ell_i}$
and $\ell_i=|q_{2i-1,2i}|$ is the order of $q_{2i-1,2i}.$  Obviously,
$$\euZ=F[t_1^{\pm 1},\ldots,t_{2s}^{\pm 1},t_{2s+1}^{\pm 1},\ldots,t_n^{\pm 1}]$$
where $t_{2s+1}=x_{2s+1},\ldots, t_n=x_n$.
\end{example}

\begin{remark}\label{torus}
In the Setting of \ref{tensorproduct} assume that all $l_i \in F^\times$.
Then
\[
     E=\euZ[x_1^{\pm 1},x_3^{\pm 1},\ldots,x_{2s-1}^{\pm 1}]\subset Q_{(\euZ)}
\]
is a maximal \'etale $\euZ$--subalgebra of $Q_{(\euZ)}$, and the pair
$(Q_{(\euZ)}, E)$  gives rise to a reductive $\euZ$--group scheme $\bGL_
{Q_{(\euZ)}}$ and its maximal torus $\mathbf{S}= R_{E/\euZ}(\bG_{m,E})$.
\end{remark}

\section{Results on $\lsl_\ell(Q)$ (mostly in good characteristic)}\label{lsl_2}

\subsection{Associative and Lie algebras.} \label{Lie-ass}
For an arbitrary associative algebra $\scA$ over a base ring $R$ we denote by
$\scA\op$ the opposite algebra: $\scA\op=\scA$ as $R$--modules but the multiplication
is given by $a\,{\cdot_{\scA\op}} \,b=b\,{\cdot_\scA}\,a$. Note that for
$\Mat_\ell(\scA)$, $\ell \in \NN_+$, namely the associative $R$--algebra of $\ell
\times \ell$-matrices over $\scA$, we have
\begin{equation}
  \label{Lie-ass1} \Mat_\ell(\scA)\op \simeq \Mat_\ell (\scA\op).
\end{equation}

 The algebra $\scA$ becomes a Lie $R$--algebra, denoted $\Lie(\scA)$, with
respect to the commutator $[a,b] = ab -ba$ as multiplication. We leave it to
the reader to check that the map
\begin{equation}
  \label{Lie-ass2}
  \iota\op \co \Lie (\scA\op)\simlgr \Lie(\scA), \quad a \mapsto -a,
\end{equation}
is an isomorphism of Lie $R$-algebras.

Abiding by the tradition notation we put
\[ \Lie \big( \Mat_\ell(\scA)\big) = \gl_\ell(\scA). \]
Thus, combining \eqref{Lie-ass1} and \eqref{Lie-ass2} we obtain a Lie algebra
isomorphism, also denoted $\io\op$,
\[ \io\op \co \gl_\ell(\scA\op) \simlgr \gl_\ell(\scA), \quad x \mapsto -x. \]

\subsection{The Lie algebra $\lsl_\ell(\scA)$.} \label{def:lsl}
Let $\scA$ be a unital associative $R$--algebra, and let $\ell \in \NN$,
$\ell \ge 2$. The derived algebra of the Lie algebra $\gl_\ell(\scA)$ of
\ref{Lie-ass}, is called the {\em special linear Lie algebra
$\lsl_\ell(\scA)$}:
\[ \lsl_\ell(\scA) = [ \gl_\ell(\scA), \gl_\ell(\scA)]. \]
Whenever we will consider $\lsl_\ell(\scA)$ in the future, it will implicitly be
assumed that $\ell \ge 2$.  Obviously, the restriction of the isomorphism
$\io\op$,
\begin{equation} \label{def:lslio}
 \io\op \co \lsl_\ell(\scA\op) \simlgr \lsl_\ell(\scA), \quad x \mapsto - x
 \end{equation}
is an isomorphism of Lie $R$--algebras. We will later need the fine structure
of $\lsl_\ell(\scA)$. First, we have
\[
    \lsl_\ell(\scA) = \{x\in \gl_\ell(\scA):  \tr(x) \in [\scA,\scA]\},
\]
where the trace $\tr(x)$ of $x\in \gl_\ell(\scA)$ is defined as usual. We
denote by $E_{ij}$, $1\le i,j\le \ell$ the usual matrix units.  It is easy to
see that
\begin{equation} \label{def:lsl00} \begin{split}
   \lsl_\ell(A) &= \scL_0 \oplus \big(\textstyle \bigoplus_{1\le i \ne j \le \ell}
            \scA E_{ij}\big), \\
    \scL_0 &= [\scA,\scA]E_{11} \oplus \{ \textstyle \sum_{i=1}^\ell a_i E_{ii}: a_i\in \scA,
                  \sum_i a_i = 0 \}.
\end{split}\end{equation} In
particular, for any unital subalgebra $S$ of $R$ the Lie $S$-algebra
$\lsl_\ell(\scA)_{(S)}$ contains
\[ \lsl_\ell(S) = \{ x\in \gl_\ell(S) : \tr(x) = 0 \}
\]
as subalgebra. We denote \begin{equation} \label{def:lsl0}
  \h_S = \lsl_\ell(S) \cap \scL_0 = \{ \textstyle \sum_{i=1}^\ell s_i E_{ii}:
             s_i\in S, \sum_i s_i = 0 \}
\end{equation}
the diagonal subalgebra of $\lsl_\ell(S)$.

We will say that a domain $R$ has {\em good characteristic for
$\lsl_\ell(\scA)$\/} if the characteristic of the fraction field of $R$ is either $0$ or $p>3$ and
such that $p$ does not divide $\ell$. In that case, if $F$ is a subfield of
$R$, the subspace $\h_F$ is an $\ad$-diagonalizable subalgebra of
$\lsl_\ell(\scA)_{(F)}$ in the sense of Subsection \ref{def:mad}, and
\eqref{def:lsl00} is the joint eigenspace decomposition of $\h_F$. \sm

It will be useful later to have a coordinate-free approach to
$\lsl_\ell(\scA)$. Namely, we let
\begin{equation}
  \label{def:lsl1} V = V_\scA = \scA \oplus \cdots  \oplus \scA =: \scA^{\oplus \ell}
\end{equation}
be the free right $\scA$-module of rank $\ell$. We denote by $\sfB = \{e_1,
\ldots, e_\ell\}$ the {\em standard basis\/} of the $\scA$--module $V$:
\begin{equation} \label{def:st}
 e_1  =  (1,0,\ldots,0),\quad \cdots \quad , e_\ell  =  (0,0,\ldots,1),
\end{equation} so that
\begin{equation}\label{def:lsl2}
 V= \textstyle \bigoplus_{i=1}^\ell e_i \scA.
\end{equation}
We let the associative algebra $\End_\scA(V)$ of $\scA$--linear endomorphisms
of $V$ act on $V$ from the left. Representing $f\in \End_\scA(V)$ by the
matrix ${\rm Mat}_\sfB(f)$ with respect to $\sfB$ provides us with an
$R$--algebra isomorphism
\begin{equation} \label{def:lsl3}
   {\rm Mat}_\sfB \co  \End_\scA(V) \simlgr \Mat_\ell(\scA)
\end{equation}
and thus also an isomorphism of the associated Lie algebras, $\gl_\scA(V)
\simeq \gl_\ell(\scA)$, and of their derived algebras,
\[ \lsl_\scA(V) := [\gl_\scA(V), \gl_\scA(V)]  \simlgr \lsl_\ell(\scA).\]

Moreover, observe
\[ e_i \scA = E_{ii} (V) \quad \text{and} \quad
    E_{ij} \co e_j \scA \simlgr e_i\scA \text{ is an isomorphism of
     $\scA$--modules.}
\]

\begin{Lemma} \label{cc} In the Setting of\/ {\rm \ref{def:lsl}} assume $\ell
\cdot 1_R\in R^\times$ and $\scA = \euZ(\scA) \oplus [\scA,\scA]$. Then \begin{align*}
      \gl_\ell(\scA) = \euZ\big( \gl_\ell(\scA) \big) \oplus \lsl_\ell(\scA) \quad \hbox{with}
      \quad\euZ\big(\gl_\ell(\scA)\big) = \euZ(\scA) E_\ell
 \end{align*}
where $E_\ell\in \gl_\ell(\scA)$ is the $\ell \times \ell$ identity matrix.
\end{Lemma}

\begin{proof} Straightforward.
\end{proof}

The following result, determining the centroid of $\lsl_\ell(\scA)$, is
folklore.

\begin{Lemma}\label{ctd-sl} Let $R$ be a commutative ring with $\frac{1}{2} \in R$, let
$\scA$ be a unital associative $R$-algebra, and let $\scL= \lsl_\ell(\scA)$
with $\ell \ge 2$.
Then for every $z\in \euZ(\scA)$ the map $\ze_z\co \scL \to \scL$, $\ze_z\big(
(x_{ij})\big) = (z x_{ij})$, is a centroidal transformation of $\scL$. The
map
\begin{equation}
  \label{ctd-sl1} \ze \co \euZ(\scA) \simlgr \Ctd(\scL), \quad z \mapsto \ze_z
\end{equation}
is an isomorphism of associative algebras. In particular, \begin{equation}
\label{ctd-sl2}
   \lsl_\ell(\scA) \text{ is fgc} \quad \iff \quad \scA \text{ is fgc.}
\end{equation}
\end{Lemma}

\begin{proof}
Since $z[a,b] = [za,b]$ it is clear that $\ze_z$ is an endomorphism of
$\scL$. It is also immediate that $\ze_z \in \Ctd(\scL)$ and that $\ze$ is an
injective homomorphism of associative algebras. Thus it remains to prove
surjectivity. Let $\chi \in \Ctd(\scL)$.

We first consider the case $\ell = 2$. For $a,b\in \scA$ we define
\begin{align*}
 e(a) &= \begin{pmatrix} 0 & a \\ 0 & 0 \end{pmatrix}, &
 f(b) & =\begin{pmatrix} 0 & 0 \\ b & 0 \end{pmatrix}, \\
 H(a,b) &= [e(a), f(b) ] =\begin{pmatrix} ab & 0 \\ 0 & -ba \end{pmatrix}, &
   h &= H(1,1) = \begin{pmatrix}  1 & 0 \\ 0 & -1\end{pmatrix}.
\end{align*}
Then
\begin{align*}
  \{ l\in \scL : [h,l] = 2 l \} &= e(\scA), \quad
  \{ l \in \scL : [h,l] = -2l \} = f(\scA),  \\
    \{ l \in \scL : [h,l] = 0 \} &= \Span \{ H(a,b) : a,b\in \scA\} =:\scL_0.
\end{align*}
It follows that $\chi$ leaves $e(\scA)$, $f(\scA)$ and $\scL_0$ invariant. In
particular we can define $\chi_\pm \in \End(\scA)$ by
\[
  \chi\big(e(a)\big) = e\big( \chi_+(a)\big), \quad
  \chi\big( f(b)\big) = f\big( \chi_-(b)\big).
\]
Denoting by $\{a\, b\, c\} = abc + cba$ the Jordan triple product of $\scA$,
we have the following multiplication rules of $\scL$: \[
    [H(a,b), e(c)] = e\big( \{a\ b\ c\}\big), \quad
    [H(a,b), f(c)] = f\big( -\{b\ a\ c\}\big).\]
They imply $\chi_\pm( \{a\ b\ c\}) = \{a\ b\ \chi_\pm(c)\}$. Define $z_\pm
\in \scA$ by $\chi_\pm(1_\scA) =  z_\pm$. Then, specializing $a=c=1_\scA$ we
obtain
\[
  2 \chi_\pm(b) = \chi_\pm\big( \{1_\scA\ 1_\scA\ b\}\big)
   = b \chi_\pm(1_\scA) + \chi_\pm(1_\scA) b = 2 z_\pm \circ b,
  \]
where $x \circ y = \frac{1}{2}(xy+ yx)$ is the Jordan algebra product of
$\scA$. Thus $\chi_\pm(b) = z_\pm \circ b$. From $H(a,b) = [e(a), f(b)]$ we
now get
 \[
      H(z_\pm \circ a, b) = [\chi(e(a)), \, f(b)] = [e(a), \chi(f(b))] =
           H(a,z_- \circ b).
 \]
 Comparing the $(11)$-entry of $H$, this proves
\[z_+ ab + a
z_+ b = a z_- b + a b z_-.\] In particular, for $a=1_\scA$ we obtain $2z_+b =
z_-b + b z_-$. Specializing $b=1_\scA$ in the formula above, this shows $z_+
= z_- =:z$, whence $2zb=zb+bz$, or $zb=bz$ for all $b\in \scA$, i.e., $z\in
\euZ(\scA)$. Finally, $\chi\big(H(a,b)\big) = H(za, b) = zH(a,b)$ proves
$\chi = \ze_z$.

Let now $\ell \ge 3$, and $i\ne j$. From
\[\scA E_{ij} = \{ l \in \lsl_\ell(\scA);
[E_{ii} - E_{jj}, \, l ] = 2 l\}\] we get $\chi(\scA E_{ij}) =\scA E_{ij}$,
allowing us to define $\chi_{ij}\in \End \scA$ by $\chi(a E_{ij}) =
\chi_{ij}(a) E_{ij}$. For distinct $i,j,p$ and $a,b\in \scA$ we have the
multiplication formula
\[
     ab E_{ij} = [[[ aE_{ij}, \, E_{jp}], \, E_{pi}], \, bE_{ij}]
\]
which implies $\chi_{ij}(ab) = \chi_{ij}(a) b = a \chi_{ij}(b)$, i.e.,
$\chi_{ij} \in \Ctd(\scA) \simeq \euZ(\scA)$ by \eqref{cbc2}. Thus, there
exists $z_{ij} \in \euZ(\scA)$ such that $\chi_{ij}(a) = z_{ij} a$. From
$[aE_{ij}, E_{jp}] = aE_{ip}$ we now obtain $z_{ij}= z_{ip}$ and from
$[E_{pi}, aE_{ij}] = a E_{pi}$ we get $z_{ij} = z_{pj}$. Hence the $z_{ij}$
are independent of $(ij)$, say $z_{ij} =:z \in \euZ(\scA)$. Finally $\chi =
\ze_z$ follows.
\end{proof}

\begin{Lemma}\label{exte}
Let $\ell$, $\ell'\in \NN$ with $\ell, \ell'\ge 2$, and let $R$ be a
commutative base ring for which $2 \cdot 1_R$, $\ell \cdot 1_R$ and
$\ell'\cdot 1_R$ are invertible in $R$. Furthermore, let $\scA$ and $\scA'$
be unital associative $R$--algebras satisfying
\[ \scA = \euZ(\scA) \oplus [\scA, \scA]\quad \hbox{and} \quad
     \scA' = \euZ(\scA') \oplus [\scA', \scA'].
\]
\sm

{\rm (a)} Let $f \co \lsl_\ell(\scA) \to \lsl_{\ell'} (\scA')$ be an $R$--linear
isomorphism of Lie algebras. The $R$--linear algebra isomorphisms $\Ctd(f)$,
$\ze_\scA$ and $\ze_{\scA'}$ of\/ {\rm \eqref{cbc0}} and\/ {\rm
\eqref{ctd-sl1}} allow us to define an isomorphism $f_\euZ \co \euZ(\scA) \to
\euZ(\scA')$ by requiring commutativity of the diagram \begin{equation}
\label{exte2} \vcenter{
\xymatrix@C=50pt{
    \euZ(\scA) \ar@{-->}[r]^{f_\euZ}  \ar[d]_{\ze_\scA}
            & \euZ(\scA') \ar[d]^{\ze_{\scA'}} \\
     \Ctd\big( \lsl_\ell(\scA)\big) \ar[r]^{\Ctd(f)^{-1}} & \Ctd\big(
     \lsl_{\ell'}(\scA') \big)
}}\end{equation}  For $z\in \euZ(\scA)$ and $X\in \lsl_\ell(\scA)$ define
\[ f_\gl \co \gl_\ell(\scA) \to \gl_{\ell'}(\scA'), \quad
     z E_\ell + X \mapsto f_\euZ(z)E_{\ell'} + f(X)
\]
Then $f_\gl$ is an isomorphism of Lie algebras.
If $\lsl_\ell(\scA)$ and $\gl_\ell(\scA)$ are viewed as $\euZ(\scA')$-algebras via the construction of {\rm \eqref{cbc0}}, then both $f$ and $f_\gl$ are $\euZ(\scA')$-linear. \sm

{\rm (b)} Let $\vphi \co \Mat_\ell(\scA) \to \Mat_{\ell'}(\scA')$ be an
isomorphism of associative algebras. The induced Lie algebra isomorphism $\vphi_\lsl \co
\lsl_\ell(\scA) \to \lsl_{\ell'}(\scA')$ obtained from $\vphi$ satisfies  $(\vphi_{\lsl})_\gl = \vphi$.
\end{Lemma}

\begin{proof}
(a) is immediate from Lemma~\ref{cc} and Lemma~\ref{ctd-sl}. For (b) we use that $\vphi$ maps the
centre $\euZ(\scA)E_\ell$ of $\Mat_\ell(\scA)$ into the centre
$\euZ(\scA')E_{\ell'}$ of $\Mat_{\ell'}(\scA')$, hence induces an $R$--linear
isomorphism $\psi \co \euZ(\scA) \to \euZ(\scA')$ by $\vphi(zE_\ell) =
\psi(z) E_{\ell'}$. Our claim is $(\vphi_\lsl)_\euZ = \psi$.

Denoting by $L$ and $L'$ the left multiplication of the associative algebras
$\Mat_{\ell}(\scA)$ and $\Mat_{\ell'}(\scA')$ respectively, we have
\[ \vphi \circ L_{z E_\ell} \circ \vphi^{-1} = L'_{\vphi(zE_{\ell})} =
L'_{\psi(z) E_{\ell'}}\] for all $z\in \euZ(\scA)$. Note that $L_{zE_\ell}$
stabilizes $\lsl_\ell(\scA)$. We denote by $(L_{zE_\ell})_\lsl$ the
restriction of $L_{zE_\ell}$ to $\lsl_\ell(\scA)$. Then $\ze_{\scA}(z) = (L_{z
E_\ell})_\lsl$. Using analogous notation for $\Mat_{\ell'}(\scA')$ and taking
the $\lsl$--components of the displayed equation above, we get
\[
 \big(  \Ctd(\vphi_\lsl)^{-1} \circ \ze_\scA \big)(z)=
    \vphi_\lsl \circ \ze_{\scA}(z) \circ \vphi_\lsl^{-1}
   = \vphi_\lsl \circ (L_{zE_\ell})_\lsl \circ \vphi_{\lsl}^{-1}
   = (L'_{\psi(z) E_{\ell'}})_\lsl = \ze_{\scA'}\big(\psi(z)\big)
\]
which proves our claim.
\end{proof}

\begin{remark}\label{exte-rem} We will later apply this Lemma in a situation
where we are given a Lie algebra isomorphism $f\co \lsl_\ell(\scA) \to
\lsl_{\ell'}(\scA')$ and \begin{enumerate}[(a)]

  \item either $f$ extends to an isomorphism $\widehat f \co \Mat_{\ell}(\scA) \to
  \Mat_{\ell'}(\scA')$ of associative algebras,

 \item or $f \circ\io\op \co \lsl_{\ell}(\scA\op) \to \lsl_{\ell'}(\scA')$
 extends to an $R$-linear isomorphism $\widehat{f \circ \io\op} \co \Mat_{\ell}(\scA\op) \to
    \Mat_{\ell'}(\scA')$ of associative algebras.
\end{enumerate}
In the first case, the Lie algebra isomorphism $f_\gl\co \gl_\ell(\scA) \to
\gl_{\ell'}(\scA')$ of Lemma~\ref{exte} is in fact an isomorphism of
associative algebras, namely $\widehat f = f_\gl$ (as maps), while in the
second case we have $\widehat{f \circ \io\op}  = (f\circ \io\op)_\gl$.
\end{remark}

\begin{defn}[{\bf AD and MAD subalgebras}] \label{def:mad} We now come to the central
concept of this paper. Let $F$ be a field. Following \cite[\S6]{CGP} we call
an $F$-subalgebra $\h$ of a Lie algebra $\scL$ over $F$ an {\em AD
subalgebra\/} if the adjoint action of each element $x\in \h$ on $\scL$ is
$F$-diagonalizable, i.e., $\scL$ admits an $F$-basis consisting of
eigenvectors of $\ad_L(x)$ for all $ x \in \h$.

A maximal AD subalgebra of $\scL$, i.e., one
which is not properly included in any other AD subalgebra of $\scL$ is called
a {\em MAD subalgebra\/}, or a {\em MAD} for short.

It is not difficult to show, see for example \cite[Lemma~8.1]{hum}, that an
AD subalgebra is necessarily abelian. Hence, AD can be thought of as an abbreviation for ``abelian
$k$-diagonalizable'' or ``$\ad$ $k$-diagonalizable''.

Let $Q$ be a quantum torus over $F$. Denote by $\euZ\big(\gl_\ell(Q)\big)$
the centre of the Lie algebra $\gl_\ell(Q).$  Assuming that $F$ is of good
characteristic for  $\gl_\ell(Q)$, by Lemma \ref{cc} we have
\[ \gl_\ell(Q) = \euZ\big(\gl_\ell(Q)\big) \oplus \lsl_\ell(Q), \quad
     \euZ\big( \gl_\ell(Q)\big) = \euZ(Q) E_\ell.\]
It follows that $\h \subset \lsl_\ell(Q)$ is an AD or a MAD of $\lsl_\ell(Q)$
if and only if $\euZ(Q)E_\ell  \oplus \h$ is an AD or a MAD of $\gl_\ell(Q)$
respectively.\sm

We next give a first example of a MAD.
\end{defn}
\begin{prop}
\label{standard-mad} Let $Q$ be a quantum torus over a field $F$ of good
characteristic for $\lsl_\ell(Q)$, $\ell \ge 2$. The subalgebra $\h_F$ of
{\rm \eqref{def:lsl0}} is a MAD of the Lie algebra\/ $\lsl_\ell(Q)_{(F)}$.
\end{prop}

We note that in case $\lsl_\ell(Q)$ is an fgc Lie torus over a field of
characteristic $0$, cf.\ \ref{qua-pro}\eqref{qua-pro-e} and \eqref{ctd-sl2},
the lemma has been proven in \cite[Cor.~5.5]{Al}. The methods of \cite{Al}
cannot be applied in our case.

\begin{proof}
It is clear that $\h=\h_F$ is an AD: the joint eigenspaces of $\h$ are the
subspaces $QE_{ij}$ and $\scL_0$ of \eqref{def:lsl00}. To show maximality,
let $d\in \lsl_\ell(Q)$ be an $\ad$ $F$-diagonalizable element commuting with
$\h$. It follows that $d\in \lsl_\ell(Q)_0 = \{ l\in \lsl_\ell(Q): [h,l]=0
\hbox{ for all }h\in \h\}$. Thus, $d=\diag(d_1, \ldots, d_\ell)$ is a
diagonal matrix. For fixed $i$, $1 \le i \le \ell$, and $q\in [Q,Q]$ we have
$qE_{ii}\in \lsl_\ell(Q)_0$ and $[d,qE_{ii}] = [d_i, q]E_{ii}$. Because
$Q=\euZ(Q) \oplus [Q,Q]$ it follows that $\ad d_i\in \End_F(Q)$ is
$F$-diagonalizable. Hence, by Lemma~\ref{diag-q}\eqref{diag-q-b}, all $d_i
\in \euZ(Q)$. Consider now the action of $d$ on an off-diagonal space
$QE_{ij}$, $i\ne j$. Clearly, $\ad d$ leaves $QE_{ij}$ invariant and acts on
$qE_{ij}$ by $(d_i - d_j) q E_{ij}$, Thus, the left multiplication by $d_i -
d_j$ is diagonalizable, forcing $d_i - d_j \in F$ by
Lemma~\ref{diag-q}\eqref{diag-q-a}. Now consider the equation
\[ [Q,Q] \ni \textstyle \sum_i d_i = (d_1 - d_2) + 2(d_2 - d_3) +\cdots
    + (n-1)(d_{\ell - 1} - d_\ell) + \ell d_\ell \in \euZ(Q)
\]
It follows that $\sum_i d_i = 0$ and that $d_\ell \in F$. Analogously, all
$d_i \in F$, and $d\in \h$ follows.
\end{proof}

\subsection{Complete orthogonal systems.}\label{como} Let $\scB$ be a unital
associative $R$--algebra. A {\em complete orthogonal system\/} (of
idempotents) in $\scB$ is a family $\scO=(e_1, \ldots, e_m)$ of elements $e_i
\in \scB$ satisfying
\begin{equation} \label{como1}
   e_i e_j = \de_{ij} e_i \quad\hbox{for $1\le i,j\le m$ and}\quad
   e_1 + \cdots + e_m = 1_\scB.
\end{equation}
In $\scB = \Mat_{\ell}(\scA)$, $\scA$ unital associative, the {\em standard
orthogonal system\/} is $\scO_{\rm st} = (E_{11}, E_{22}, \ldots, \allowbreak
E_{\ell \ell})$, where the $E_{ij}$ are the usual standard matrix units. But
also $(E_{11} + E_{22}, E_{33}, \ldots, E_{\ell \ell})$ is a complete
orthogonal system. Part (a) of the following Lemma~\ref{comoc} says that
there is a natural bijection between complete orthogonal systems in
$\Mat_\ell(\scA)$ and decompositions of $V_\scA$ as a direct sum of
$\scA$--modules.

\begin{lemma} \label{comoc}
Let $\scA$ be a unital associative $R$--algebra, and let $V=V_\scA=
\scA^{\oplus \ell}$ be the right $\scA$--module of \eqref{def:lsl1}. We
identify $\scB = \Mat_\ell(\scA) \equiv \End_\scA(V)$ using
\eqref{def:lsl3}.\ms

\begin{inparaenum}[\rm (a)]

\item \label{comoc-a}  Let $\scO=(e_1, \ldots, e_m)$ be a complete orthogonal
    system in $\scB$. Define $V_i = e_i(V)$, $1\le i \le m$. Then $V$ decomposes as
    $V=V_1 \oplus \cdots \oplus V_m$ where each $V_i$ is a right
    $\scA$--module. Conversely, let $V=V_1 \oplus \cdots \oplus V_m$ be a decomposition of $V$ as
a direct sum of $\scA$--modules. Define $e_i \in \scB$ as the canonical
projection of $V$ onto $V_i\subset V$. Then $(e_1, \ldots, e_m)$ is a
complete orthogonal system in $\Mat_\ell(\scA)$.

The constructions $\scO \rightsquigarrow V=V_1 \oplus \cdots \oplus V_m$ and
$V=V_1 \oplus \cdots \oplus V_m \rightsquigarrow \scO$ defined above are
inverses of each other. \ms

\item \label{comoc-b}  Let $\scO=(e_1, \ldots, e_m)$ and $\scO'=(e_1',
    \ldots, e'_{m'})$ be complete orthogonal systems in $\Mat_\ell(\scA)$, inducing
    by \eqref{comoc-a} decompositions $V_\scA = V_1 \oplus \cdots \oplus V_m = V_1'
    \oplus \cdots \oplus V'_{m'}$. Let
\[
          \scD(\scO) = \{ f\in \End_\scA (V_\scA) : f(V_i) \subset V_i , 1 \le i \le m\}
              = \textstyle \bigoplus_{i=1}^m \End_\scA(V_i)
\]
and define $\scD(\scO')$ analogously. Suppose that all $\End_\scA(V_i)$ and
$\End_\scA(V'_j)$ are connected. Then the following are equivalent for $g\in
\scB^\times$.

\begin{inparaenum}[\rm (i)]
   \item \label{comoc-bi} $g \scD(\scO) g^{-1} =
       \scD(\scO')$.

   \item \label{comoc-bii} $m=m'$ and there exists a permutation $\si\in \SS_m$ such that
       $ge_i g^{-1} = e'_{\si(i)}$, $1 \le i \le m$.

   \item \label{comoc-biii} $m=m'$ and there exists a permutation $\si \in \SS_m$ such that
   $g(V_i) = V'_{\si(i)}$ for $1\le  i \le m$.
\end{inparaenum}
\sm

\item \label{comoc-c} Let $\scA = Q$ be a quantum torus. Let $\scO_{\rm st} =
    (E_{11}, \ldots, E_{\ell \ell})$ be the standard orthogonal system in
    $\scB = \Mat_\ell(Q)$, and let $\scO=(e_1, \ldots e_\ell)$ be another
    complete orthogonal system in $\scB$ with associated decomposition $V_Q =
    V_1 \oplus \cdots \oplus V_\ell$, $V_i = e_i(V)$ for $1 \le i \le \ell$. Define $\scD(\scO_{\rm
    st})$ and $\scD(\scO)$ as in \eqref{comoc-b}. Then the following are
    equivalent.

 \begin{inparaenum}[\rm (i)]
   \item \label{comoc-ci} Each $V_i$, $ 1\le i \le \ell$,  is a cyclic $Q$--module.

   \item \label{comoc-cii} Each $V_i$, $1 \le i \le \ell$, is a free $Q$--module of rank $1$.

    \item \label{comoc-ciii} There exists $g\in \scB^\times$ such that $g \scD(\scO_{\rm st}) g^{-1} =
    \scD(\scO)$ and each $\End_Q(V_i)$ is connected.
 \end{inparaenum}

\noindent Assuming \eqref{comoc-ciii} holds, let $g\in \scB^\times$ be as in
\eqref{comoc-ciii} and let
$$\h_{\rm st} = \{\sum_{i=1}^\ell s_i E_{ii}: s_i
\in F, \sum_i s_i=0\}$$ be the standard MAD of $\lsl_\ell(Q)$ and define $\h
\subset \scD(\scO)$ analogously. Then the automorphism $\Int(g)$ of
$\lsl_\ell(Q)$ maps $\h_{\rm st}$ onto $\h$. In particular, $\h$ is also a
MAD of $\lsl_\ell(Q)$. \
\end{inparaenum}
\end{lemma}

\begin{proof}
\eqref{comoc-a} That $V=\sum_{i=1}^m V_i$ follows from the second equation in
 \eqref{como1}, and that the sum is direct from the first. The converse is
 equally straightforward.

\eqref{comoc-b} Assume \eqref{comoc-bi}. Since $(ge_i g^{-1})_{1\le i \le m}$
is a complete orthogonal system of $\scB$ contained in $\scD(\scO')$, it
follows from our connectedness assumption that each $e_i\in \scO$ is a sum of
some of the $e'_i \in \scO'$, and that distinct idempotents in $\scO'$ are
used for each $e_i$. Hence $m \le m'$. By symmetry, $m'\le m$, whence $m=m'$.
The remaining part of \eqref{comoc-bii} is now clear. The implications
$\eqref{comoc-bii} \implies \eqref{comoc-biii} \implies \eqref{comoc-bi}$ are
immediate from the definitions.

\eqref{comoc-c}. A cyclic $Q$--module is free since $Q$ is a domain. Thus
\eqref{comoc-ci} $\iff$ \eqref{comoc-cii}. For the proof of \eqref{comoc-cii}
$\implies$ \eqref{comoc-ciii} put $V_{i, {\rm st}}= E_{ii}(V)$. We know
$V_{i, {\rm st}} \simeq Q$ by \eqref{def:lsl2}. Also, by assumption, there
exist $Q$--linear isomorphisms $g_i \co V_{i, {\rm st}} \to V_i$, $1\le i \le
\ell$. Hence $g= g_1 \oplus \cdots \oplus g_\ell$ is an invertible
endomorphism of $V$ such that $g(V_{i, {\rm st}}) = V_i$, $1\le i \le \ell$.
It then follows from \eqref{comoc-b} that $g \scD(\scO_{\rm st}) g^{-1} =
\scD(\scO)$ (note that \eqref{comoc-b} can be applied since $Q\simeq
\End_Q(V_i)$ is connected by \ref{qua-pro}\eqref{qua-pro4}. The implication
\eqref{comoc-ciii} $\implies$ \eqref{comoc-cii} also follows from
\eqref{comoc-b}.
\end{proof}


\section{Isomorphisms between two Lie algebras of type ${\rm A}$ over  rings}

\subsection{Torsion bijections}\label{sec:torbi} We start by reviewing some
of the techniques of non-abelian \v Cech cohomology used later on. \sm

Let $\bG$ be a smooth affine group scheme over a (commutative, unital) ring
$R$. The pointed set of non-abelian \v Cech cohomology on the
 \'etale site of ${X=\rm Spec}(R)$ with coefficients in $\bG$, is denoted by
 $H_{\et}^1(X, \bG)$.  This pointed  set measures the isomorphism classes of torsors over $X$ under
 $G$ (see  \cite[ Ch.\,IV \S1]{milne} and \cite{DG} for basic definitions and references). Abusing notation a bit we
 will identify the set of isomorphism classes of $\bG$-torsors over $X$ with $H^1_{\et}(R,\bG)$.
\sm

Recall that any morphism  $\bG\to \bH$ of group schemes induces a natural map
$H^1_{\et}(R,\bG) \to H^1_{\et}(R,\bH)$. If $[E]\in H^1_{\et}(R,\bG)$ we will
denote its image in $H^1_{\et}(R,\bH)$ by  $[E_\bH]$. \sm

For a $\bG$-torsor $E$   we denote by  ${^E\bG}$ the twisted form of $\bG$ by $E.$ This is
a smooth affine group scheme over $X$. Recall that
according to \cite[III.2.6.3.1]{Gi} there exists a natural bijection
$$
\tau_E: H_{\et}^1(X, {^E\bG}) \to H_{\et}^1(X, \bG),$$ called the {\em
torsion bijection\/}, which takes the class of the trivial torsor under
$^E\bG$ to the class of $E$. \sm

Let $[E]\in H^1_{\et}(R,\bG)$. Any exact sequence
\[
    1\longrightarrow \bG \stackrel{\psi}{\longrightarrow} \bH
              \longrightarrow \bF \longrightarrow 1 \] of
smooth affine group schemes  induces a commutative diagram
$$
\begin{CD}
 \bF(R) @>>> H^1_{\et}(R,  {^E\bG}) @>{\psi_E}>> H^1_{\et}(R, {^{E_\bH}\bH}) @>>>
            H^1_{\et}(R,\bF)\\
@V{\wr}VV @V{\tau_E}VV @V{\tau_{E_H}}VV @VVV \\
^E\bF(R)=\bF(R) @>>> H^1_{\et}(R,\bG) @>{\psi}>> H^1_{\et}(R,\bH) @>>>
            H^1_{\et}(R,\bF). \\
\end{CD}
$$

\begin{lemma}\label{chase} Using the notation of\/ {\rm \ref{sec:torbi}}, the torsion bijection $\tau_E$ induces a bijection between ${\rm
Ker}(\psi_E)$ and the fiber $\psi^{-1}(\psi(E))=\psi^{-1}(E_H)$.
\end{lemma}
\begin{proof} This follows from an easy diagram chase.
\end{proof}

\subsection{Azumaya algebras} \label{Azu}
Let $\scA$ be an Azumaya algebra over a (commutative, associative, unital) ring
$R$.
If  $\scA$ has rank $\ell^2$,   it is a
twisted form of the matrix algebra $\Mat_\ell(R)$. Since
\[
   \bAut_R(\Mat_\ell(R))\simeq\bPGL_{\ell,R}
\]
(see \cite[Chapter IV, Proposition 2.3]{milne}), the elements of the pointed
set $H_{\et}^1(R,\bPGL_{\ell,R})$ are in one-to-one correspondence with the
isomorphism classes of Azumaya algebras over $R$ of degree $\ell$ (the bijection
is given by twisting).
It follows that  $\scA\simeq {^{\xi}\Mat_\ell(R)}$ for some class $[\xi]\in
H^1_{\et}(R, \bPGL_{\ell,R})$ and that  $ \bAut_R(\scA)\simeq \bPGL_{\scA}$.

\subsection{Automorphisms of $\lsl_{\scA}$.} \label{aut-pre} (a) Let $\scA$ be an Azumaya
algebra over $R$. Every $\vphi \in \bAut({\scA})(R)$ leaves the Lie algebra
\[    \lsl_{\scA} := [\scA,\scA] \]
invariant, thus induces an automorphism $\vphi_\lsl$ of
$\lsl_{\scA}$. Since the construction $\vphi \mapsto \vphi_\lsl$ is functorial, it
gives rise to a homomorphism
\begin{equation} \label{aut-pre1}
 \bPGL_{\scA} \to \bAut(\lsl_{\scA}) \end{equation}
which is injective since $\lsl_{\scA}$ generates $\scA$ as an associative algebra. \sm

(b) Assume $\scA$ has an anti-automorphism $\ka$. Then \begin{equation*}
   \ka_\lsl \co \lsl_{\scA} \to \lsl_{\scA}, \quad x \mapsto - \ka(x)
\end{equation*}
is an automorphism of $\lsl_{\scA}$. Again by functoriality of the construction,
this gives rise to an element of $\bAut(\lsl_{\scA})$, also denoted $\ka_\lsl$. In
case $\scA = \Mat_\ell(R)$ we use the transpose as anti-automorphism and put
\begin{equation}\label{aut-pre2}
    \ta \co \lsl_\ell(R) \to \lsl_\ell (R), \quad x \mapsto -\,{^t x}
\end{equation}
As before, this gives rise to an automorphism in $\bAut(\lsl_\ell(R))$, also
denoted $\ta$. \sm

(c) Putting together the maps in (a) and (b) we have constructed
homomorphisms of $R$--group schemes
\begin{equation}\label{aut-pre3}
   \bPGL_{2,R} \to \bAut\big(\lsl_2(R)\big)
\end{equation}
and for $m\ge 3$
\begin{equation}
  \label{aut-pre4}
      \bPGL_{m,R} \rtimes (\ZZ/2\ZZ)_R \to \bAut\big( \lsl_m(R)\big); \quad
          (\vphi, \eps) \mapsto \vphi_\lsl \circ \ta^\eps
\end{equation}
The reader will easily check that both maps are injective homomorphisms (one
needs the assumption $m\ge 3$ to get injectivity in the second case).

\begin{theorem}\label{aut-th} Assume $R$ is a domain of good characteristic for
$\lsl_m$ containing a field $F$. Then the maps \eqref{aut-pre3} and \eqref{aut-pre4} are
isomorphisms of $R$--groups schemes.
\end{theorem}

\begin{proof}
Since all groups involved are obtained from $F$ by base change, and since
base change preserves isomorphism, we may assume that $R=F$ is a field of
good characteristic. In this case it is shown in \cite[4.7]{St} (or see
\cite[Thm.~IX.5]{jake:lie} for $F$ of characteristic $0$ and
\cite[p.~67]{sel} for characteristic $>0$) that for any field extension $E/F$ the maps \eqref{aut-pre3} and
\eqref{aut-pre4}  evaluated at the $E$--points are isomorphisms of abstract
groups. In particular this holds for the algebraic closure of $F$. A standard
fact in the theory of group schemes (see \cite[Proposition $22.5$]{KMRT}) 
then proves the result.
\end{proof}
\ms
\subsection{Consequences}\label{con}

We will derive some consequences of Theorem~\ref{aut-th}.
Let again $\scA$ be an Azumaya algebra over $R$.  Assume that $R$ is a
domain. Since then $\scA$ has constant rank as an $R$--module, it is a twisted
form of $M_\ell(R)$ for some $\ell$.

 For the remainder of this section we
will assume that $R$ contains a field of good characteristic for
$\lsl_\ell(R)$. We consider the corresponding Lie $R$-algebra
$\lsl_{\scA}=[\scA,\scA]$.
By a standard twisting argument, it follows from Theorem~\ref{aut-th} that we
have an $R$-group scheme isomorphism
\[\bPGL_{\scA}\simeq \bAut_R(\lsl_{\scA}) \]
if $\ell = 2$, and  the exact sequence  of $R$--group schemes.
\[
 1\longrightarrow \bPGL_{\scA}\longrightarrow \bAut_R(\lsl_{\scA})
    \longrightarrow (\mathbb{Z}/2\mathbb{Z})_R\longrightarrow 1
\]
if $\ell \geq 3.$

Evaluating at $R$-points we have
\[
1\longrightarrow {\rm PGL}_{\scA}=\bPGL_{\scA}(R)\hookrightarrow
{\rm Aut}_R(\lsl_{\scA})=\bAut_R(\lsl_{\scA})(R)\longrightarrow
\mathbb{Z}/2\mathbb{Z}.
\]
Note that this last map is trivial when $\ell = 2.$ We will say that an automorphism $\phi\in {\rm Aut}_R(\lsl_{\scA})$ is {\em
inner\/} if it is in the image of ${\rm PGL}_{\scA}$. Otherwise we say that $\phi$
is {\em outer}.

\begin{theorem}\label{Vcon} 
Let $\phi\in {\rm Aut}_R(\lsl_{\scA})$.
Then the following holds:
\begin{enumerate}[\rm (a)]
\item\label{Vcon-a}  If $\phi$ is inner, it is the restriction of a unique
    automorphism $\psi:\scA\to \scA$.

\item\label{Vcon-b} If $\phi$ is outer, it is the restriction of the negative of a unique
    anti-automorphism $\psi:\scA\to \scA$.
\end{enumerate}
\end{theorem}

\begin{proof} The subset $Y=\bAut_R(\lsl_{\scA}) \setminus  \bPGL_{\scA}$ is a closed subscheme of
$\bAut_R(\lsl_{\scA})$ consisting of outer automorphisms of $\lsl_{\scA}$.
The group $\bPGL_{\scA}$ acts simply transitively on $Y$
by left multiplication. Thus $Y$ is a $\bPGL_{\scA}$-torsor. It is trivial if and only if $\lsl_{\scA}$
has at least one outer automorphism.

Along the same lines, let $X$ be the scheme of anti-automorphisms of ${\scA}$:
\[
X=\{ \,\psi:  \scA\to \scA\ | \,\  \psi \,\, \text{\rm is bijective and} \,  \ \psi(xy)=\psi(y)\psi(x) \mbox{\ for all \ } x,y\in \scA\,\}.
\]
We observe that if $\psi$ is an anti-automorphism of $\scA$ and $\ell\geq 3$ then $-\psi|_{\lsl_{\scA}}$
is an outer automorphism of $\lsl_{\scA}$.
The automorphism group $\bPGL_{\scA}$ of ${\scA}$ acts simply transitively on $X$ on
the left (because the action is simply transitively in the split case). Thus,
$X$ is also a $\bPGL_{\scA}$-torsor. As before, if $\ell\geq 3$ then $X$ is a trivial torsor, i.e.,
there exists at least one anti-automorphism of $\scA$, if and only if $\lsl_{\scA}$
has at least one outer automorphism.

Note that the natural restriction map $\lambda: X\to Y$ which takes $\psi$
into $-\psi|_{\lsl_{\scA}}$ is an isomorphism of torsors (because it is an
isomorphism in the split case).

\eqref{Vcon-a} Let $\phi$ be inner. Then there exists a unique
$g\in\bPGL_{\scA}(R)$ whose image in $\bAut_R(\lsl_{\scA})$ is $\phi$. This element $g$
corresponds to the automorphism of $\scA$ which we  denote by $\psi:\scA\to \scA$. By
construction, the restriction $\psi|_{\lsl_{\scA}}$ is $\phi$. The uniqueness of
such a $\psi$ is immediate, since $\lsl_{\scA}$ generates $\scA$ as
associative
algebra.

\eqref{Vcon-b} Let now $\phi$ be outer. Then since $\lambda$ is bijective
there exists a unique $\psi\in X(R)$ such that $\lambda(\psi)=-\psi|_{\lsl_{\scA}}$
equals $\phi$.
\end{proof}

\begin{corollary}\label{order}
Let $\scA$ be an Azumaya algebra over $R$, thus a twisted form of
$\Mat_\ell(R)$. Assume $\ell \ge 3$. Then $\scA \simeq \scA\op$ if and only
if $\lsl_\scA$ has an outer automorphism.
\end{corollary}

\begin{proof}
Assume $\scA \simeq \scA\op$. We have seen in \ref{aut-pre}(b) that an
isomorphism $\ka \co \scA \to \scA\op$ induces an automorphism $\ka_\lsl$ of
$\lsl_\scA$. We claim that $\ka_\lsl$ is an outer automorphism. It suffices
to show this after some base change. Take any base change $R\to S$ which
splits $\scA$, i.e., $\scA_S \simeq \Mat_\ell(S)$. Then $\ka$ is a composite
of an inner conjugation, say $int(g)$, and the transpose $\ta$. It is
well-known that for $\ell \ge 3$ the restriction of $-\ta$ to $\lsl_\ell$ is
an outer automorphism (because it induces the automorphism $-1$ of 
the corresponding root system), while the restriction of $int(g)$ to $\lsl_\scA$ is an
inner automorphism. Therefore $\ka_\lsl$ is indeed outer.

The other direction follows from Theorem~\ref{Vcon}\eqref{Vcon-b}.
\end{proof}

\begin{theorem}\label{AorAop}
Let $\scA,\scA'$ be Azumaya algebras over our domain $R$  of the
same rank $l^2$. If\/ $\lsl_{\scA}\simeq \lsl_{\scA'}$ as $R$--algebras, then
\smallskip

\begin{inparaenum}[\rm (a)]
\item \label{AorAop-a} $\scA'\simeq \scA$ if $\ell=2$;

\item \label{AorAop-b} $\scA'\simeq \scA$ or $\scA'\simeq \scA\op$ if $\ell \geq 3$.
    \end{inparaenum}\end{theorem}

We note that the converses in (a) and (b) are  obvious.

\begin{proof} We abbreviate $H^1= H^1_{\et}$ and let $[\xi],[\xi']\in  H^1(R,\bPGL_{\ell,R})$
be the classes  corresponding to $\scA$ and $\scA'$. Because $\scA={^\xi \Mat}_\ell(R)$
the Lie algebra $\lsl_{\scA}$ is the twisted form of $\lsl_\ell(R)$ by $\xi$. \ms

{\it Case $\ell=2$}: Recall from Theorem~\ref{aut-th} that \[
   \bAut_R(\lsl_2(R)) = \bPGL_{2,R},\]
i.e., the Lie algebra $\lsl_2(R)$ and the Azumaya algebra $\Mat_2(R)$ have
the same automorphism group scheme. Therefore, if the twisted forms $\lsl_{\scA}$
and $\lsl_{\scA'}$ of the Lie algebra $\lsl_2(R)$ are isomorphic, then
$[\xi]=[\xi']$ implies $\scA\simeq \scA'$ as $R$--algebras. \ms

{\it Case $\ell \geq 3$}: Recall from Theorem~\ref{aut-th} that
 \[
    \bAut_R( \lsl_\ell(R)) = \bAut_R (\bPGL_{\ell,R}) = \bPGL_{\ell,R}\rtimes (\mathbb{Z}/2\mathbb{Z})_R.
\]
Hence we
get the exact sequence of group schemes
\begin{equation}\label{two}
 1\longrightarrow \bPGL_{\ell,R}\longrightarrow \bAut_R(\lsl_\ell(R))
    \longrightarrow (\mathbb{Z}/2\mathbb{Z})_R\longrightarrow 1
\end{equation}
which induces a canonical map
 \[
H^1(R,\bPGL_{\ell,R})\stackrel{\psi}{\longrightarrow}
H^1(R,\bAut_R(\lsl_\ell(R))=H^1(R,\bAut_R(\bPGL_{\ell,R}))
\]
in cohomology. By assumption, $\psi([\xi])=\psi([\xi'])$
(because $\lsl_{\scA} \stackrel{R}{\simeq} \lsl_{\scA'}$). If $[\xi]=[\xi']$ then the
twisted Azumaya algebras $\scA={^{\xi}\Mat_\ell(R)}$ and $\scA'={^{\xi'}\Mat_\ell(R)}$
are isomorphic as $R$--algebras, and we are done.

Assume therefore $[\xi]\not=[\xi']$. Twisting (\ref{two}) by $\xi$ we get an
exact sequence of group schemes
$$
 1\longrightarrow \bPGL_{\scA}\longrightarrow \bAut_R(\bPGL_{\scA})
\longrightarrow (\mathbb{Z}/2\mathbb{Z})_R\longrightarrow 1
$$
which in turn induces an exact sequence
\begin{equation}\label{automorphisms}
\bPGL_{\scA}(R)\hookrightarrow \bAut_R(\bPGL_{\scA})(R)\longrightarrow
\mathbb{Z}/2\mathbb{Z}\longrightarrow H^1(R,\bPGL_{\scA})
\stackrel{\phi}{\longrightarrow}
H^1(R,\bAut_R(\bPGL_{\scA}))
\end{equation}
of pointed sets. According to  Lemma~\ref{chase} there is a natural
one-to-one correspondence between the fiber $\psi^{-1}(\psi(\xi))$ and the
kernel of $\phi$. Since $[\xi] \ne [\xi']\in \psi^{-1}(\psi(\xi))$, the class
$[\xi']$ corresponds to some nontrivial element in ${\rm Ker}(\phi)$. On the
other hand,  by \eqref{automorphisms}, ${\rm Ker}(\phi)$ consists of at most
two elements implying  $|{\rm Ker}(\phi)|=2$. We then get from
(\ref{automorphisms}) that $\bPGL_{\scA}(R)= \bAut_R(\bPGL_{\scA})(R)$, i.e., every
automorphism of $\PGL_{\scA}$ over $R$ is inner. By
Corollary~\ref{order}, this implies ${\scA}\not\simeq {\scA}\op$. \sm

Now we observe that the opposite Azumaya algebra $\scA\op$ also corresponds to
some element in $H^1(R,\bPGL_{\ell,R})$, call it $[\xi\op]$. We have
$[\xi]\not=[\xi\op]$ because $\scA\not\simeq \scA\op$. Furthermore, it is
straightforward to verify that the map $\bGL_{\scA} \to \bGL_{\scA\op}$, given by
$x\to x^{-1}$, is an isomorphism of group schemes which induces in turn an
isomorphism $\bPGL_{\scA}\to \bPGL_{\scA\op}$. This implies that
$\psi([\xi])=\psi([\xi\op])$, hence the class $[\xi\op]$ also corresponds to
the nontrivial element ($\not=[\xi]$) in the fiber $\psi^{-1}(\psi([\xi]))$.
Since $ |\psi^{-1}(\psi([\xi]))|=| {\rm Ker}(\phi)|=2$, necessarily
$[\xi']=[\xi\op]$ in $H^1(R,\bPGL_{\ell,R})$, implying $\scA'\simeq \scA\op$.
\end{proof}

As a by-product of the proof we obtain that the converse statement in
Corollary~\ref{order}
also holds.

\begin{corollary}\label{inner}
Let $\scA$ be an Azumaya algebra over $R$ of rank $\ell^2 \ge 4$, where $R$ has
good characteristic for $\lsl_\ell(R)$. Suppose $\scA\not\simeq \scA\op$.
Then every automorphism of the Lie $R$--algebra $\lsl_{\scA}$ is inner. If in
addition $\Pic(R)=1$ then every automorphism of $\lsl_{\scA}$
is the
restriction of the conjugation map by an element $a\in \scA^{\times}$.
\end{corollary}

\begin{proof} Without loss of generality we may assume that $\ell\geq 3$.
We use the notation of Theorem~\ref{AorAop}.
Since $\scA\not\simeq \scA\op$ we have $[\xi]\not=[\xi\op]$.
Therefore, in sequence (\ref{automorphisms}) the kernel of
the canonical map $\phi$ consists of two elements and this implies
$\bPGL_{\scA}(R)=\bAut_R(\bPGL_{\scA})(R)=\Aut(\scA)$.
Thus every automorphism of $\lsl_{\scA}$ is induced by an automorphism of $\scA$.
It remains to note that if $\Pic(R)=1$ the natural map $\scA^{\times}\to
\bPGL_{\scA}(R)$ is surjective.
\end{proof}

In the following sections we apply all these results to the case of fgc
quantum tori viewed over their centres.


\section{$\lsl_\ell(Q)$ for $Q$ an fgc quantum torus}\label{slfgc}
Our analysis will use the following description of isomorphisms between
special linear Lie algebras over rings. The theorem below is
an immediate consequence of our Theorem~\ref{Vcon} and is
originally  due to Jacobson-Seligman
when the Lie algebras in question are defined over fields.

\begin{theorem} \label{JSTheorem} Let $\scA$ and
$\scA'$ be Azumaya algebras of the same rank over a domain $R$ containing a field
of good characteristic.
Let $f \co \lsl_{\scA'} \to \lsl_{\scA}$ be an
$R$--linear isomorphism of Lie algebras. Then
 after possibly replacing $\scA$ by $\scA\op$, the map $f$ uniquely extends
 to an $R$--linear isomorphism $\tilde f \co \scA' \to \scA$
 of associative algebras.
\end{theorem}
\begin{proof}
Since it is clear that $f \circ \io\op$ is a
Lie  $R$--algebra isomorphism, we only need to show that $f$  extends to a
map $\tilde f \co \scA' \to \scA$ which is either an
isomorphism or the negative of an anti-isomorphism of associative algebras.
By Theorem~\ref{AorAop},  after possibly replacing $\scA$ by $\scA\op$ we can assume that
there exists an isomorphism
$g\co \scA' \to \scA$ of associative
$R$--algebras. Since the composition $f\circ g^{-1}|_{\lsl_{\scA}}$ is an $R$-automorphism of
$\lsl_{\scA}$, by Theorem~\ref {Vcon} it can be extended uniquely to either an automorphism
or the negative of anti-automorphism  of $\scA$
 and the assertion follows.
\end{proof}
\begin{remark}
For the special case of central-simple algebras over a field
$K$ of characteristic $0$, the
theorem is proven in \cite[Thm. IX.5]{jake:lie}. The characteristic
$0$ case can easily be extended to
positive characteristic by using the description of ${\rm Aut}\, \lsl_\ell(K)$ in \cite[p.~67]{sel}.
\end{remark}

We will apply this result in case $\scA$ and $\scA'$ are matrix algebras over
fgc quantum tori.

\begin{theorem} \label{jaseco} Let $\scL = \lsl_\ell(Q)$ and
$\scL'=\lsl_{\ell'}(Q')$ where $Q$ and $Q'$ are fgc quantum tori over a field
$F$  of good characteristic for $\scL$ and $\scL'$. We denote by
$\euZ=\euZ(Q)$ the centre of $Q$, and let $f \co \scL' \to \scL$ be an
$F$--linear isomorphism of Lie algebras. Then the following hold. \sm

\begin{inparaenum}[\rm (a)]

 \item \label{jaseco-c} After possibly replacing $Q$ by $Q\op$, the map $f$ uniquely extends
 to a $\euZ$--linear isomorphism $\tilde f \co \Mat_{\ell'}(Q')_{({\euZ})} \to \Mat_\ell(Q)_{(\euZ)}$
 of associative algebras.

  \item\label{jaseco-a} $\ell' = \ell$.

\end{inparaenum}
\end{theorem}

\begin{proof}
(a)
By \ref{ctd-sl} the centre $\euZ$ of $Q$ is isomorphic to the
centroid of $\scL$ -- we will take this as an identification. Applying
\ref{cbc}, Example 1, our map $f$ is a $\euZ$--linear isomorphism
\[
    f \co \scL'_{(Z)} = \lsl_{\ell'}(Q'_{(\euZ)}) \simlgr
        \scL_{(\euZ)} = \lsl_{\ell} (Q_{(\euZ)}).
\]
Observe that the two Azumaya algebras $\Mat_\ell(Q)$ and $\Mat_{\ell'}(Q')$
over $\euZ$ have the same rank.
Indeed, passing to a proper \'etale base extension we may assume that they are split.
Let $m^2$ and $(m')^2$ be their ranks. Then ${\rm dim}\,\scL=m^2-1$ and ${\rm dim}\,\scL'=
(m')^2-1$ implying $m=m'$. Now
Theorem~\ref{JSTheorem} applied to
$\scA=\Mat_\ell(Q)_{(\euZ)}$ and $\scA'=\Mat_{\ell'}(Q')_{(\euZ)}$
implies the claim.

(b)
Let $K$ be the fraction field of $\euZ$ (recall from
\ref{qua-pro}\eqref{qua-pro2} that $\euZ$ is a Laurent polynomial ring over
$F$).
Then after possibly replacing $Q$ by $Q\op$ we can assume that
$\tilde{f}$ uniquely extends to a $K$-linear isomorphism
$$\tilde{f}_K:M_{\ell'}(Q')\otimes_{\scZ} K \simlgr  M_\ell(Q)\otimes_{\scZ} K$$
of associative algebras.
But $D'=
Q'_{(\euZ)} \ot_\euZ K$ and $D = Q_{(\euZ)} \ot_\euZ K$ are central division
algebras over $K$.
Now 
by the theorem of Wedderburn (see \cite[Theorem $2.1.6$]{He}) 
we have $D'\simeq D$ and
$\ell'=\ell$. 
\end{proof}

The spirit of  Theorem~\ref{jaseco}, in a sense that will be made precise,  will allow us to reduce questions of isomorphisms of
Lie algebras to questions of isomorphisms of associative algebras. These are handled in the following result.

\begin{theorem}\label{isomQQ'}
Let $Q$ and $Q'$ be fgc quantum tori over the field $F$ and let $\scA
=\Mat_{\ell}(Q)$ and $\scA'=\Mat_{\ell}(Q')$. We assume that the
characteristic of $F$ is {\rm very good}: it is either $0$, or $p>3$ where
$p$ does not divide $\ell$ or any of the $q_{ij}$, $q'_{ij}$ occurring in the
quantum matrices of $Q$ and $Q'$ respectively with respect to a canonical
presentation. We let $\euZ$ be the centre of $Q$, and assume that $f \co
\scA' \to \scA$ is an $F$--linear isomorphism of associative algebras. As in
{\rm \ref{cbc}} we will view $f$ as a $\euZ$--linear isomorphism
$\scA'_{(\euZ)} \to \scA_{(\euZ)}$. The isomorphism $f$ induces an
isomorphism $\phi\co G'\to G$ of the associated reductive group schemes
$G'=\bGL_{A'_{(\euZ)}}$ and $G=\bGL_{A_{(\euZ)}}$ over $\euZ$.
Let $T'$ (resp. $T$) be the split torus in $G'$
    (resp.\ $G$) corresponding to the diagonal matrices in $A'_{(\euZ)}={\rm
    M}_{\ell}(Q'_{(\euZ)})$ (resp. in $A_{(\euZ)}= {\rm M}_{\ell}(Q_{(\euZ)})$).
Then $\phi(T')$ and $T$ are conjugate in $G$.
\end{theorem}

\begin{proof}
Clearly, $T'$ and $T$ are
generically maximal split tori in the sense of \cite{CGP}, i.e., $T'_{K}$ and
$T_{K}$ are maximal split tori in $G_{K}$ and $G'_{K}$ where $K=\tilde \euZ$
is the fraction field of $\euZ$ (because $\widetilde Q=
Q_{(\euZ)}\otimes_{\euZ} K$ and $\widetilde{Q'}_{(\euZ)}\otimes_{\euZ} K$ are
central division algebras over $K$ by \ref{qua-pro}\eqref{qua-pro5}). Note
that the centralizer of $T$ (resp.\ $T'$) in $G$ (resp.\ in $G'$) is
isomorphic to the reductive group scheme $C_G(T)\simeq
\bGL_{Q_{(\euZ)}}\times \ldots \times \bGL_{Q_{(\euZ)}}$ over $\euZ$ (resp.
$C_{G'}(T')\simeq \bGL_{Q'_{(\euZ)}}\times\ldots\times \bGL_{Q'_{(\euZ)}}$).

Since $T'$ and $T$ are generically maximal split tori, the proof of
\cite[Proposition 8.1]{CGP} shows that an obstacle for conjugacy of
$\phi(T')$ and $T$ in $G$ is given by a class
$$
[\xi]\in H^1_{Zar}(\euZ,C_{G}(T))=H^1_{Zar}(\euZ,\bGL_{Q_{(\euZ)}})
\times\ldots\times
H_{Zar}^1(\euZ,
\bGL_{Q_{(\euZ)}})
$$
and that
$$
^{\xi}(C_{G}(T))\simeq C_{G}(\phi(T'))\simeq C_{G'}(T')\simeq \bGL_{Q'_{(\euZ)}}
\times\ldots\times \bGL_{Q'_{(\euZ)}}.
$$
Let $\xi=(\xi_1,\ldots,\xi_\ell)$. Since
$^{\xi_i}\bGL_{Q_{(\euZ)}}\simeq \bGL_{Q'_{(\euZ)}}$,
for the proof of
conjugacy of $\phi(T')$ and $T$
 it suffices to show that $\xi_i=1$. \sm

We now recall some general facts from the theory of reductive group schemes.
Let $H$ be any reductive group scheme over $\euZ$, $S\subset H$  a maximal
torus and  $[\xi]\in H^1_{\et}(\euZ,H)$. Let $N=N_H(S)$ and $W=N/S$.
\smallskip

\begin{enumerate}[(i)]

\item {\it By \cite[Lemma $8.2$]{CGP}, $W$ is a finite \'etale $\euZ$-group and 
the canonical map $H^1_{\et}(\euZ,W)\to H^1_{\et}(K,W)$ is injective; in particular $H^1_{Zar}(\euZ,W)=1$.}
 
\item \label{isomQQ'-i} {\it According to \cite[Remarque $3.2.5$]{Gi} if the twisted group scheme $^{\xi}H$ contains
    a  maximal torus then the class $[\xi]$ is in the image of
    $H^1_{\et}(\euZ,N_H(S))\to H^1_{\et}(\euZ,H)$.}
\smallskip

\item \label{isomQQ'-ii} 
 {\it Let $S'\subset H$ be another maximal torus. Then the transporter ${\rm {\bf Trans}}_H(S,S')$ 
 is an $N$-torsor over $\euZ$, hence it corresponds to a unique class $\lambda  \in H^1_{\et}(\euZ,N)$. If $S,S'$
 are conjugate over $K$ then $\lambda$ is in the image
 of $H^1_{\et}(\euZ,S)\to H^1_{\et}(\euZ,N)$. If in addition $H^1(\euZ,S)=1$ then
$S$ and $S'$ are conjugate over $\euZ$.}
\end{enumerate}
\smallskip

For the proof of the first statement in (iii) we refer to \cite[Lemma $8.3$]{CGP}. The second assertion follows from (i) and the commutativity of the following diagram:
$$
\begin{CD}
 H^1_{\et}(\euZ, S) @>>> H^1_{\et}(\euZ, N) @>>>
            H^1_{\et}(\euZ,W)\\
@VVV @VVV @VVV \\
H^1_{\et}(K,S) @>>> H^1_{\et}(K,N) @>>>
            H^1_{\et}(K,W). \\
\end{CD}
$$

We will apply these facts in the following two  situations. Let $m$ be an integer such that 
the matrix algebra ${\rm M}_m(\euZ)$ is the split form of $Q_{(\euZ)}$. Take the pair 
$(E,S)$ consisting of a maximal commutative \'etale subalgebra $E\subset Q_{(\euZ)}$
and the corresponding torus $S=R_{E/\euZ}({\rm {\bf G}}_{m,E})$ constructed 
in Remark~\ref{torus}. Fix any embedding $E\hookrightarrow {\rm M}_m(\euZ)$ (for instance,
we can take the regular representation of $E$). This also gives rise to a closed embedding
$S\hookrightarrow H=\bGL_{m,\euZ}$. Note that the torus $S$ determines uniquely
the subalgebra $E$. Indeed, $E$ is the centralizer of $S(\euZ)$ in ${\rm M}_m(\euZ)$.
When we view $E$ as a subalgebra in $Q_{(\euZ)}$ we will denote it by $E_Q$. Similarly the 
corresponding torus
$S$ will be denoted by $S_Q$. 
Consider now another embedding $E\to {\rm M}_m(\euZ)$. We denote its image by $E'$
and the corresponding torus by $S'$.

\begin{lemma} 
The two \'etale subalgebras $E,E'\subset {\rm M}_m(\euZ)$ are conjugate under the action of
$\bGL_m(\euZ)$.
\end{lemma} 
\begin{proof} It suffices to show the conjugacy of $S$ and $S'$. The obstacle for conjugacy of 
$S,S'$ is the $N$-torsor ${\rm {\bf Trans}}_H(S,S')$ where $N=N_{\bGL_{m,\euZ}}(S)$.
By the Skolem-Noether Theorem (see~\cite[Theorem $1.4$]{KMRT}) the algebras $E$ and $E'$
are conjugate over $K$. Furthermore, we have 
$$H^1_{\et}(\euZ,S)=H^1_{\et}(E,{\rm{\bf G}}_{m,E})={\rm Pic}(E)=1,$$ 
because $E$ is a Laurent polynomial ring  and hence has trivial Picard group.
The claim now follows from (\ref{isomQQ'-ii}).
\end{proof}

We next pass to adjoint groups.
\begin{lemma} 
We keep the above notation. Let $\overline{H}=\bPGL_{m,\euZ}$ and let $\overline{S}$ be the image
of $S$ under the canonical map $\bGL_{m,\euZ}\to \bPGL_{m,\euZ}$.
There exists a class $[\theta]\in H^1(\euZ,\overline{S})$ such that
$^{\theta}\bPGL_{m,\euZ}\simeq \bPGL_{Q_{(\euZ)}}$. 
\end{lemma}
\begin{proof} Since $\bPGL_{Q_{(\euZ)}}$ has a maximal torus, by  
(\ref{isomQQ'-i}) the group scheme $\bPGL_{Q_{(\euZ)}}$ is a twisted form
of $\bPGL_{m,\euZ}$ by some cocycle $\lambda$ with coefficients in
$\overline{N}=N_{\bPGL_{m,\euZ}}(\overline{S})$. The cohomological class 
$[\lambda]$ corresponds to the maximal torus $\overline{S}_Q\subset \bPGL_{Q_{(\euZ)}}$
where $\overline{S}_Q$ is the image of $S_Q\subset \bGL_{Q_{(\euZ)}}$ under the 
canonical map $\bGL_{Q_{(\euZ)}}\to\bPGL_{Q_{(\euZ)}}$.
The argument in
\cite[Theorem $3.1$]{Ch} shows that there is another  closed embedding
$\overline{S}\hookrightarrow \bPGL_{m,\euZ}$, whose image will
be denoted by $\overline{S}'$, such that $[\lambda]$ is equivalent
to some class $[\lambda']\in H^1(\euZ,\overline{S}')$. Let $S'\subset \bGL_{m,\euZ}$
be the preimage of $\overline{S}'$. Its centralizer $E'$ in $Q_{(\euZ)}$
is a maximal commutative \'etale subalgebra of $Q_{(\euZ)}$ isomorphic to $E$. Since by the above lemma
$E$ and $E'$ are conjugate over $\euZ$,   so are the maximal tori
$S,S'$. This in turn implies  the conjugacy of $\overline{S}$ and $\overline{S}'$ in
$\bPGL_{m,\euZ}$. Thus the cocycle $\lambda'$ is equivalent to
some cocycle $\theta$ with coefficients in $\overline{S}$.  
\end{proof} 

Assume for a moment that $E$ also admits an embedding into $Q'_{(\euZ)}$.
Then as above the group scheme $\bPGL_{Q'_{(\euZ)}}$ is a twisted from of $\bPGL_{m,\euZ}$
with some cocycle $\theta'$ with coefficients in $\overline{S}$. Consider the exact
sequence 
$$
H^1(\euZ,S)=1\longrightarrow H^1(\euZ,\overline{S})\longrightarrow
H^2(\euZ,\bG_{m,\euZ})={\rm Br}(\euZ).
$$
The images of $[\theta]$ and $[\theta']$ in ${\rm Br}(\euZ)$
coincide, because the base change morphism ${\rm Br}(\euZ)\hookrightarrow 
{\rm Br}(K)$ is injective and
$Q_{(\euZ)}$ and $Q'_{(\euZ)}$ are isomorphic over $K$. It follows
that $[\theta]=[\theta']$ and this implies $Q_{(\euZ)}\simeq Q'_{(\euZ)}$;
in particular $\xi_i=1$ as required.

To sum up: to finish the proof of Theorem~\ref{isomQQ'}
what is left  to show
is that $E$ admits an embedding in $Q'_{(\euZ)}$, i.e., there exist elements
$y_1,y_3,\ldots,y_{2s-1}\in (Q'_{(\euZ)})^{\times}$ which commute and such
that
$$
y_1^{\ell_1}=t_1,\ y_3^{\ell_3}=t_3,\ldots,y_{2s-1}^{\ell_s}=t_{2s-1}.
$$
Here the positive integers $\ell_1,\ldots,\ell_s$ and variables $t_1,\ldots,t_n$
are the same as in Example~\ref{tensorproduct} applied to $Q$.

Fix a presentation $Q'=\bigoplus_{\la \in \La' } k x'^\la$. As usual it induces
a  grading of $Q'$.
The isomorphism $f\co \scA'_{(\euZ)}\to \scA_{(\euZ)}$ induces an isomorphism
\[
     \Mat_\ell(Q'_{(\euZ)})\otimes_{\euZ} K=\Mat_\ell(Q'_{(\euZ)}\otimes_{\euZ}
           K)\to \Mat_\ell(Q_{(\euZ)})\otimes_{\euZ}
             K=\Mat_\ell(Q_{(\euZ)}\otimes_{\euZ} K).
\]
From the theory of
central simple algebras over fields we know that the last implies
that $Q'_{(\euZ)}\otimes_{\euZ} K\simeq
Q_{(\euZ)}\otimes_{\euZ} K$. Hence, there exist commuting elements
$z_1,z_3,\ldots,z_{2s-1}\in Q'_{(\euZ)}\otimes_{\euZ} K$ such that
$z_{2i-1}^{l_i}=t_{2i-1}$ for all $i=1,\ldots,s$, where we view $t_i \in
Q'_{(\euZ)}$. Choose an element $r\in \euZ$ such that $rz_{2i-1}\in
Q'_{(\euZ)}$ for all $i$. Then
\begin{equation}\label{highest}
(rz_{2i-1})^{l_i}=r^{l_i}t_{2i-1}.
\end{equation}
Recall that $Q'_{(\euZ)}$ has a natural grading transferred from $Q'.$
\ms
Let $u_{2i-1}$ (resp. $a$) be the highest
homogeneous component of $rz_{2i-1}$ (resp. $r$) with respect to any order on $\Lambda$. Taking the highest
components on the left and on the right  in (\ref{highest}) we get
$(u_{2i-1})^{l_i}=a^{l_i}t_{2i-1}$. Note that $a$ is some monomial in the
centre $\euZ^{\times}$ of $Q'_{(\euZ)}$, hence it commutes with $u_{2i-1}$.
Then the elements
$$y_{2i-1}=a^{-1}u_{2i-1}\in (Q'_{(\euZ)})^{\times},\  i=1,3,\ldots,2s-1,$$ have the required properties.
\end{proof}

\begin{corollary}\label{cor-conj}
Let $\scL = \lsl_\ell(Q)$ and $\scL'=\lsl_{\ell'}(Q')$ where $Q$ and $Q'$ are
fgc quantum tori over a field of very good characteristic, let
$\mathfrak{m}\subset \lsl_\ell(Q)$ (resp. $\mathfrak{m}'\subset
\lsl_{\ell'}(Q')$) be the MAD consisting of diagonal matrices with entries in
$F$, and let $f \co \scL' \to \scL$ be an $F$--linear isomorphisms of Lie
algebras. Then $f(\mathfrak{m}')$ is conjugate to $\mathfrak{m}$.
\end{corollary}

\begin{proof}
  This follows by combining Theorem~\ref{jaseco} and Theorem~\ref{isomQQ'}.
\end{proof}

\section{Specialization of quantum tori}\label{specialization}

Our main method of dealing with non-fgc quantum tori is {\em
specialization\/} which we develop in this section. In this section $k$
denotes a field of characteristic $0$.

\begin{prop}   \label{speca}
Let $\scF$ be a finite subset of a field $k$ of characteristic $0$ consisting
of non-zero elements and let $\ell \in \NN_+$. Then there exists a finitely
generated subring $R \subset k$ and a maximal ideal $\mmm \triangleleft R$
such that
\begin{enumerate}[\rm (a)]
  \item\label{speca-a} $\scF \subset R\setminus \mmm$, and

  \item \label{spec-b} $R/\mmm$ is a finite field of characteristic $p >
      \ell$.
\end{enumerate} \end{prop}

\begin{proof} We can assume that $f^{-1} \in \scF$ for every $f\in \scF$.
Let $R = \ZZ[\scF] \subset k$ (resp.\ $C=\bbQ[\scF] \subset k$) be the
$\ZZ$--subalgebra (resp.\ $\bbQ$--subalgebra) generated by $\scF$. By the
Noether Normalization Lemma there exist algebraically independent $u_1,
\ldots, u_s \in C$ over $\bbQ$ such that $C$ is integral over its subring
$\bbQ[u_1, \ldots, u_s]$ and of finite type as a $\bbQ[u_1, \ldots,
u_s]$--module, say with generators $c_1, \ldots, c_t$. \sm

{\em Observation.} For finitely many polynomials $q_1, \ldots, q_m \in
\bbQ[u_1, \ldots, u_s]$ there exists $n\in \NN_+$ such that $q_1, \ldots, q_m
\in \ZZ[\frac{1}{n}][u_1, \ldots, u_s]$ (here $\ZZ[\frac{1}{n}]$ is the
localization of $\ZZ$ in $\{ n^a: a \in \NN\}$). \sm

Apply this observation to the coefficients of the minimal polynomials of the
integral elements $c_1, \ldots, c_t$ and to the coefficients appearing in the
linear combinations expressing the elements of $\scF$ in $\sum_{i=1}^t
\bbQ[u_1, \ldots, u_s]c_i$. This yields that there exists $n\in \NN_+$ such
that \begin{enumerate}[(1)]

\item\label{speca-1} all coefficients of the minimal polynomials of $c_1,
    \ldots,  c_t$ belong to $E:=\ZZ[\frac{1}{n}][u_1, \ldots, u_s]$ and

\item \label{speca-2} $R:=\ZZ[\scF]\subset E[c_1, \ldots, c_t] =: D \subset C
    =
    \bbQ[\scF]$.
\end{enumerate}
Because of \eqref{speca-1}, each $c_i$ is integral over $E$, whence $D/E$ is
an integral extension of finite type. Now choose a prime number $p$ such that
$p \not| n$ and $p > \ell$. The ideal $\pp\triangleleft E$, generated by $p$
and the $u_1, \ldots, u_s$, has the property that $E/\pp \simeq
\ZZ[\frac{1}{n}]/\langle p\rangle$, where $\langle p \rangle = p
\ZZ[\frac{1}{n}]\subset \ZZ[\frac{1}{n}]$ is the ideal generated by $p$,
whence $\langle p \rangle \subset \ZZ[\frac{1}{n}]$ and therefore also $\pp
\subset E$ are maximal ideals.
Since $D$ is  integral over $E$, there exists a maximal ideal $\nn
\triangleleft D$ lying over $\pp \triangleleft E$. By construction, $D/\nn$
is a finite field of characteristic $p$.
Recall $R \subset D$, and put $\mmm = R \cap \nn$. Then $R/\mmm
\hookrightarrow D/\nn$. So $R/\mmm$ is a finite subring of the field $D/\nn$,
whence a field itself. It remains to observe that $f\not\in \mmm$ for every
$f\in \scF$ because $f$ is a unit in $D$.  \end{proof}

\begin{corollary} \label{fgc-spec}
Let $Q$ be a quantum torus over a field $k$ of characteristic $0$, let
$q=(q_{ij}) \in \Mat_n(k)$ be the quantum matrix associated with a
coordinatization of $Q$ and write $Q=\bigoplus_{\la \in \La} k x^\la$ as in
\eqref{qua-pro11}. Further, let $a_1, \ldots, a_t \in k\setminus \{0\}$,
$b_1, \ldots, b_m \in Q$ be non-zero elements and let $g_1, \ldots , g_p \in
\Mat_\ell(Q)$ be non-zero matrices. \sm

Then there exists a finitely generated subring $R< k$ and a maximal ideal
$\mmm \ideal R$ with the following properties:
\begin{enumerate}[\rm (i)]
\item all $a_i$ and $q_{ij} \in R$, all $b_1, \ldots, b_m $ lie in the
    unital graded subalgebra $\scA= \bigoplus_{\la \in \La} R x^\la$ of $Q$,
    and all $g_1, \ldots, g_p \in \Mat_\ell(\scA)$.

\item Denoting by $\overline{\phantom{0}}$ the canonical quotient map, we
    have
  \begin{enumerate}[\rm (a)]
      \item $\bar R = R/ \mmm$ is a finite field of very good
          characteristic $p>0$ for $\lsl_\ell(\bar{\scA})$,

    \item all $\bar a_i$ and $\bar q_{ij} \in \bar R$ are non-zero, hence roots of unity, 
    
    \item ${\bar {\scA}} = \bigoplus_{\la \in \La} \bar R x^\la$ is an fgc
        quantum torus over $\bar R$ with quantum matrix $(\bar q_{ij})$
        and $p \not | \, [{\bar \scA} : \euZ({\bar \scA})]$;

    \item all $\bar b_i \ne 0$ in $\bar \scA$;

    \item all $\bar g_i$ are non-zero in $\Mat_\ell(\bar \scA)$.
  \end{enumerate}

\end{enumerate}
\end{corollary}

\begin{proof} Let $\scF$ be the finite subset of $k$ consisting of all $a_i$, $q_{ij}$,
all non-zero $k$-coefficients of   $b_1,\ldots,b_m$ and $g_1,\ldots,g_p$ and
their inverses (the $k$-coefficients are taken with respect to a natural
$k$-basis of $Q$ and $\Mat_{\ell}(Q)$). Let $R$ and $\mmm$ be as in
Proposition~\ref{speca}. Then all claims follow immediately from $\scF
\subset R \setminus \mmm$. \end{proof}

\begin{lemma}\label{fin-generat}
Let $Q$ be a quantum torus over a field $k$ of characteristic $0$. Then the
Lie algebra  $\lsl_\ell(Q)$ is finitely generated over $k$.
\end{lemma}

\begin{proof}
This fact is known for all Lie tori (\cite[Thm.~5]{n:tori}). In our concrete
case, it can be proven as follows. We fix a parametrization $Q =
\bigoplus_{\al \in \ZZ^n} k x^\la$, in particular this gives us coordinates
$x_i$, $1 \le i \le n$ corresponding to the standard basis of $\ZZ^n$. It is
straightforward to check that $\lsl_\ell(Q)$ is generated by $\{ E_{ij},
x_p^{\pm 1} E_{ij}: 1 \le i \ne j \le \ell\}$.
\end{proof}

We have seen in \ref{qua-pro}\eqref{qua-pro-e} that a quantum torus is fgc if
for one coordinatization, equivalently for all coordinatizations, the entries
of the associated quantum matrix are roots of unity. Hence, in a
coordinatization of a non-fgc quantum torus with quantum matrix $q=(q_{ij})$
at least one of the $q_{ij}$ is not a root of unity.

\begin{theorem}\label{SpecializingIsomorphisms}
Let $Q = \bigoplus_{\la \in \La } k x^\la$ and $Q' = \bigoplus_{\la' \in \La'
} k y^{\la'}$ be  non-fgc quantum tori over a field $k$ of characteristic $0$
with associated quantum matrices $q=(q_{ij})\in \Mat_n(k)$ and $q'=(q'_{ij})
\in \Mat_{n'}(k)$. We assume that we are given:

\quad $\bullet$ non-zero elements $b_1, \ldots, b_m \in Q$ and non-zero
elements   $b'_1, \ldots, b'_m \in Q'$;

\quad $\bullet$ non-zero elements $g_1, \ldots , g_s \in \gl_\ell(Q)$, and
non-zero elements $g'_1, \ldots , g'_{s'} \in \gl_{\ell'}(Q')$;

\quad $\bullet$ a $k$-linear isomorphism $f: \lsl_\ell(Q) \to
\lsl_{\ell'}(Q')$ of
    Lie algebras.
\ms

\noindent Then there exists a subring $R< k$ and a maximal ideal
    $\mmm \ideal R$ with the following properties:
\begin{enumerate}[\rm (i)]
\item all $q_{ij} \in R$ and all $q'_{ij}\in R$,

\item all $b_1, \ldots, b_m $ lie in the unital graded subalgebra $\scA=
    \bigoplus_{\la \in \La} R x^\la$ of $Q$,
    and all $g_1, \ldots, g_s $ are in
    $\gl_\ell(\scA)$,

\item all $b'_1,\ldots, b'_m$ lie in the unital graded subalgebra $\scA'=
    \bigoplus_{\la' \in \La'} R y^{\la'}$ of $Q'$, and all
    $g'_1,\ldots,g'_{s'}\in \gl_{\ell'}(\scA')$.

\item  $f(\lsl_{\ell}(\scA))=\lsl_{\ell'}(\scA')$.

   \item Denoting by $\overline{\phantom{0}}$ the canonical quotient map,
       we have
  \begin{enumerate}[\rm (a)]
      \item $\bar R = R/ \mmm$ is a finite field of very good
          characteristic for $\lsl_\ell(\bar{\scA})$ and
          $\lsl_{\ell'}(\bar{\scA})$,

    \item all $\bar q_{ij}, \bar q'_{ij} \in \bar R$ are non-zero, hence
        roots of unity, 

    \item $\bar \scA = \scA/\mmm \scA = \bigoplus_{\la \in \La} \bar R
        x^\la$ and $ {\bar \scA}' = \bigoplus_{\la' \in \La'} \bar R
        y^{\la'}$  are fgc quantum
        tori over $\bar R$ with associated quantum matrices $(\bar q_{ij})$ and $(\bar q_{ij}')$, 

     \item $\bar f: \lsl_{\ell}(\bar \scA)\to \lsl_{\ell'}(\bar \scA')$ is an
         $\bar R$-isomorphism of the corresponding Lie $\bar
         R$--algebras,

\item all ${\bar b_i}$  are non-zero in $\bar \scA$ and all ${\bar b'_i}$
    are non-zero in $\bar \scA'$,

    \item all ${\bar g_i}$ are non-zero in $\gl_\ell(\bar \scA)$  and all
        ${\bar g'_i}$ are non-zero in $\gl_{\ell'}(\bar \scA')$.
  \end{enumerate}
\end{enumerate}
\end{theorem}

\begin{proof} By definition, the coordinates of $0 \ne b \in Q$ are the non-zero
$s_\la \in k$ when $b$ is written as $b= \sum_{\la \in \La} s_\al x^\la$. The
coordinates of $0\ne g=\sum_{i,j} g_{ij} E_{ij}\in \Mat_\ell(Q)$ are the
coordinates of the non-zero $g_{ij} \in Q$. The coordinates of $0\ne b'\in
Q'$ and $0\ne g'\in \Mat_{\ell'}(Q')$ are defined analogously. We choose
finite generating systems $S\subset \lsl_\ell(Q)$ and $S'\subset
\lsl_{\ell'}(Q')$ as in Lemma~\ref{fin-generat}, and put $g= f^{-1}$. We let
$\scF\subset k$ consist of
\begin{enumerate}[-]

\item all $q_{ij}$, $q_{ij}'$, the coordinates of all $b_i$, $b_i'$, $g_i$,
    $g_i'$ together with their inverses,

 \item  the elements $(q_{ij} - 1)^{-1}$ and $(q'_{ij} -1)^{-1}$
     whenever $q_{ij}$ or $q'_{ij}$ is not a root of unity, and

 \item the coefficients of all elements in $f(S)$ and $g(S')$.
\end{enumerate}

We now apply Proposition~\ref{speca} for this $\scF$ but with the $\ell$
there replaced by $\max\{\ell, \ell'\}$. This provides us with $(R,\mmm)$ as
required in the Theorem. Indeed, by construction we have
$f(\lsl_{\ell}(\scA))\subset \lsl_{\ell'}(\scA')$ and similarly
$g(\lsl_{\ell'}(\scA'))\subset \lsl_{\ell}(\scA)$. It follows  that
$f(\lsl_{\ell}(\scA))= \lsl_{\ell'}(\scA')$ (because $g$ is the inverse for
$f$), so that (iv) holds. The remaining claims follow immediately from $\scF
\subset R \setminus \mmm$.
\end{proof}

\section{Some preliminaries for step 3 of the proof of the main theorem}
\label{sec:E}

\subsection{Setting I}\label{setI}  In this section we use the following setting:
\begin{itemize}
  \item $Q$ is a quantum torus over a field $F$ with grading group
      $\La\simeq \ZZ^n$;

   \item We fix a basis $\beps$ of $\La$; the $\beps$-trace, the
       corresponding $\ZZ$-grading, and the $\beps$-degree will all be
       taken with respect to the fixed $\beps$. We will therefore simply
       write $\deg$ instead of $\deg_\beps$. We identify $\La = \ZZ^n$ via
       $\beps$, and define $\La^+ = \NN^n\subset \ZZ^n$.

 \item Corresponding to $\beps$ there exists a coordinatization of $Q$ as $Q =
     \bigoplus_{\la \in \ZZ^n} F x^\la$. We let $0\ne x_i \in Q$ denote the element that
     corresponds to the $i$th basis vector $\eps_i$ in $\beps$, and put
     $Q^+ = \bigoplus_{\la \in \La^+} Q^\la$. Note that $Q^+$ is a unital
     subring of $Q$.

 \item For $\ell \ge 2$ we let $V$ be a free right $Q$--module of rank
     $\ell$. We fix a basis $e_1, \ldots e_n$ of $V$ so that we can write
     $V=\bigoplus_{i=1}^\ell e_iQ$. We put $V^+ = \bigoplus_{i=1}^\ell e_i
     Q^+$.

 \item We say that {\em $x_i$ divides $q\in Q^+$\/} if every $\la \in \supp
     (q)$ has the form $\la = (\la_1, \ldots,
      \la_n)$ with respect to $\beps$ and $\la_i >0$. In this case $q
      x_i^{-1} \in Q^+$.

 \item Any $0 \ne v\in V$ can be uniquely written as $v=\sum_{i=1}^\ell
     e_iq_i$ with $q_i \in Q$. We put
    \[
      \deg(v) = \max \{ \deg(q_i) : q_i \ne 0 \}. \\
    \]
    We refer to the $q_i$ as the coordinates of $v$.
\end{itemize}

\begin{lemma}\label{dpm} Let $0 \ne v\in V$ and let $0 \ne q \in Q$. Then $vq \ne 0$ and
\[
   \deg (vq) = \deg(v) + \deg(q).
\]
\end{lemma}

\begin{proof}
We write $v= \sum_{i=1}^\ell e_i q_i$ as above. The first claim is obvious:
$q_i q \ne 0 \iff q_i \ne 0$, and at least one coordinate $q_i \ne 0$. For
the non-zero $q_i$ we have $\deg(q_i q) = \deg(q_i) + \deg(q)$ by
\eqref{deg-prop1}. Hence, if for simpler notation $\deg(v) = \deg(q_1)$ then
$\deg(q_iq) \le \deg(q_1 q)$ for all $i$ with $q_i \ne 0$.\end{proof}

\subsection{Setting II}\label{setII}
We continue with the Setting of \ref{setI}. In addition, we fix a non-zero
$Q$--submodule $U \subset V$, and put $U^+ = U
     \cap V^+$. Observe $U^+ \ne \{0\}$ since $qQ \cap Q^+ \ne \{0\}$ for
     any $0\ne q \in Q$. Since $\deg(u) \ge 0$ for every $0 \ne u \in U^+$ there exist
 {\em minimal\/} vectors   $u_0 \in U^+$ such that
        \[ \deg(u_0) \le \deg(u) \quad \hbox{for all $0\ne u\in U^+$}.
         \]

 We call $u\in U^+$ {\em indivisible\/} if $u\ne 0$ and  the
      coordinates of $u$ do not have a common divisor $x_i$, $1 \le i \le
      n$, in $Q^+$. It is obvious that indivisible vectors exist. Moreover,
      \[ \hbox{\em every minimal vector $u_0 \in U^+$ is indivisible.}\]
      Indeed, write
       $u_0 = \sum_{j=1}^\ell e_j q_j$ with all $q_j \in U^+$. Assume
      that all $q_j$ are divisible by some $x_i$.
     Then $0 \ne u_0
       x_i^{-1} = \sum_{j=1}^\ell e_j q_j x_i^{-1} \in U^+$ since all $q_j
       x_i^{-1} \in Q^+$. But $\deg(u_{\min} x_i^{-1}) = \deg(u_0) + \deg(x_i^{-1})
       = \deg(u_0) - 1< \deg(u_0)$, contradiction.

\begin{lemma}\label{lem811} Assume Setting {\rm \ref{setII}}, and let $u_0 \in U^+$ be indivisible, and let $q\in Q$. Then
    \[ u_0q \in U^+ \quad \iff \quad q\in Q^+.\]
\end{lemma}

\begin{proof}
We only need to show that the left-hand side implies the right-hand side.
 We reason by contradiction. Assume $q\not\in Q^+$. We write $u_0 = \sum_{j=1}^\ell e_j
q_j$ with all $q_j \in Q^+$. We choose a minimal $\la =(\la_1, \ldots,
\la_n)\in \La^+$ such that $qx^\la \in Q^+$. Here ``minimal'' means that
$qx^\la x_i^{-1} \not\in Q^+$ for all $i$.  Since $q\not\in Q^+$, some $\la_i
>0$. Let $I$ be the ideal of the ring $Q^+$ generated by $x_i$. Note
$I=x_iQ^+ = Q^+ x_i$. By definition of an indivisible element, at least one
of the
      coordinate $q_j$ of $u_0$ is not divisible by $x_i$ in $Q^+$ (for $i$
      with $\la_i >0)$. To simplify notation, assume this is $q_1$, i.e.,
      $q_1 \not\in I$. Then the following holds.
\begin{enumerate}[(i)]
   \item $qx^\la \not\in I$. Otherwise, $qx^\la = q'x_i$ for some $q' \in
       Q^+$, whence $qx^\la x_i^{-1} \in Q^+$, and so $qx^\mu \in Q^+$ for
       $\mu = \la - \eps_i \in \La^+$.

  \item $q_1 q x^\la \in Q^+$ because $q_1 \in Q^+$ and $qx^\la \in Q^+$.
      Furthermore, $x^\la \in I$ because $\la_i >0$. Hence $q_1 q x^\la \in
      I$. Thus $x_1 | q_1 q x^\la$, but $x_1 \not | q_1$ and $x_1 \not | q
      x^\la$. This contradicts (iii) below.

  \item $Q^+/I$ is a subring of the quantum torus with associated quantum
      matrix $q'=(q_{ij})$ where ${2 \le i,j\le n}$. It is therefore a domain.
\end{enumerate}
Our assumption $q\not\in Q^+$ has thus led to a contradiction.
\end{proof}

\begin{corollary}
  \label{genco} Let $u_0 \in U^+$ be an indivisible vector. Then
  \[ U=u_0 Q \quad \iff \quad U^+ = u_0 Q^+. \]
\end{corollary}

\begin{proof}
If $U=u_0 Q$ and $u^+\in U^+$, then $u^+=u_0 q$ with $q\in Q$ and
Lemma~\ref{lem811} shows $q\in Q^+$. Conversely, if $U^+ = u_0 Q^+$ and $u\in
U$ is arbitrary, we can take $x^\la \in Q$ such that $ux^\la \in U^+$, whence
$ux^\la = u_0 q$ for some $q\in Q^+$. But then $u=u_0q (x^\la)^{-1} \in u_0
Q$. \end{proof}

\begin{lemma}\label{lem812}
  In addition to the Setting\/ {\rm \ref{setII}}  assume that $U$ admits a
  complement: $V=U \oplus U'$ for some $Q$--subspace $U'$. Let $u_0 \in
  U^+$ be a minimal vector and put $Q^{++} = \bigoplus_{\la \in \NN^n \setminus \{0\}} Q^\la$.
  Then $u_0 \not\in v(Q^{++}) $ for any $v\in V^+$.
\end{lemma}

\begin{proof}
  Assume to the contrary that $u_0 = vq $ for some $q\in Q^{++}$ and $v\in V^+$.
  Decompose $v= u + u'$ with $u \in U$ and $u'\in U'$. Then $u_0 = uq + u'q$
  shows $u'q = 0$. So without loss of generality we can assume $v\in U^+$.
   Now apply Lemma~\ref{dpm}  to get
  $$
  \deg(v) + \deg(q) = \deg(vq) =
  \deg(u_0) \le \deg(v)$$
  (because $u_0$ is minimal) and hence $\deg(q) = 0$, i.e., $q\in F^{\times} \cdot 1_F$,
  contradiction. \end{proof}

\section{Proof of the main theorem}

\subsection{Setting and plan of the proof.}\label{set1} Throughout this section, $k$ is a base field of characteristic $0$
and $Q$ and $Q'$ are non-fgc quantum tori. We assume that they are
coordinatized as
\[\textstyle
   Q = \bigoplus_{\la \in \La } k x^\la\quad \hbox{and}\quad
      Q' = \bigoplus_{\la' \in \La' } k y^{\la'}\]
for $\La = \ZZ^n$ and $\La'= \ZZ^{n'}$ with associated quantum matrices
$q=(q_{ij})\in \Mat_n(k)$ and $q'=(q'_{ij}) \in \Mat_{n'}(k)$. We assume that
\[
    f\co \lsl_{\ell'}(Q') \to \lsl_{\ell}(Q)
\]
is a $k$--linear isomorphism. We apply Lemma~\ref{exte} to extend $f$ to an
isomorphism
\[ f_\gl \co \gl_{\ell'}(Q') \to \gl_{\ell}(Q).\]
In the first step of the proof (Proposition~\ref{step1}) we will show $\ell =
\ell'$. Next, it will follow from Proposition~\ref{step2} that we may assume
that
\[
   \phi = f_\gl \co \Mat_{\ell}(Q'))\to \Mat_{\ell}(Q)
\]
is an isomorphism of associative algebras. In the final step of the proof we
will establish that if $\h'\subset \Mat_{\ell'}(Q')$ (resp. $\h\subset
\Mat_{\ell}(Q)$) is the standard MAD of $\Mat_{\ell'}(Q')$ (resp.
$\Mat_\ell(Q)$) then $\phi(\h')$ is conjugate to $\h$ by an element of
$\GL_\ell(Q)$.

\begin{proposition}\label{step1} In the setting\/ {\rm \ref{set1}} we have
$\ell = \ell'$.
\end{proposition}

\begin{proof}
According to Theorem~\ref{SpecializingIsomorphisms} there exist a subring
$R\subset k$ and a maximal ideal $\mmm\triangleleft  R$ such that $f$ induces
an $\overline{R}$-isomorphism
$$
\bar{f} \co \lsl_{\ell'}(\bar{\scA'})\to \lsl_\ell(\bar{\scA}).
$$
Since $\bar{\scA}$ and $\bar{\scA'}$ are quantum tori over $\bar{R}$
of fgc type, by Theorem~\ref{jaseco} we have $\ell=\ell'$.
\end{proof}

\begin{prop}  \label{step2} Consider  the following two maps: \begin{enumerate}[\rm (a)]
   \item $f_\gl \co \gl_{\ell}(Q') \to \gl_{\ell}(Q),$

   \item the extension of $f\circ \io\op \co \lsl_{\ell}((Q')\op) \to
       \lsl_\ell(Q)$ to a Lie algebra isomorphism
       \[ (f\circ \io\op)_\gl \co \gl_{\ell}((Q')\op) \to \gl_\ell(Q).\]
\end{enumerate} Then one of these maps is an isomorphism of the underlying associative $k$--algebras.
\end{prop}
\begin{proof} Assume the contrary. Then there exist $g_1,g_2,g_3,g_4\in\Mat_\ell(Q')$
such that
\[
   f_\gl(g_1g_2)\not=f_\gl(g_1)f_\gl(g_2) \quad \hbox{and} \quad
f_\gl(-g_3g_4)\not=f_\gl(g_4)f_\gl(g_3).
 \]
Put
$$
g_1'=f_\gl(g_1g_2),\;  g_2'=f_\gl(g_1),\;
    g_3'=f_\gl(g_2),\; g_4'=g_1'-g_2'g'_3 $$ and
similarly
\[
g_5'=f_\gl(-g_3g_4),\; g_6'=f_\gl(g_3), \; g_7'=f_\gl(g_4), \; g_8'=g_5'-g_7'g_6'.
\]
Recall that by Lemma~\ref{cc} one has the decomposition
\begin{align*}
      \gl_\ell(Q) = \euZ(Q)E_\ell  \oplus \lsl_\ell(Q)
 \end{align*}
and similarly
\begin{align*}
      \gl_\ell(Q') = \euZ(Q')E_\ell  \oplus \lsl_\ell(Q').
 \end{align*}
So every element $g_i$ (resp. $g_i'$) can be written as the sum $g_i=
q_{i}E_\ell+\tilde{g}_i$ (resp. \ $g_i'=q'_iE_{\ell}+\tilde{g}'_i$) where
$q_i$ (resp.\ $q'_i$) is in the centre of $Q$ (resp.\ $Q'$) and $\tilde{g}_i$
(resp.\ ($\tilde{g}'_i$) is a sum of commutators of elements of
$\gl_\ell(Q)$ (resp.\ $\gl_{\ell'}(Q')$). We add to our list of elements
$g_1,\ldots,g_4$ (resp. $g'_1,\ldots,g'_8$) all their components arising  in
the above two decomposition (including elements appearing in the writing of
$\tilde{g}_i,\tilde{g}'_i$ as sums of commutators).

We now apply Theorem~\ref{SpecializingIsomorphisms} with these data. This
provides us with a subring $R\subset k$ and a maximal ideal $\mmm
\triangleleft R$ satisfying the many conclusions of loc.\ cit. In particular,
denoting by \[ f_\scA \co \lsl_{\ell}(\scA') \to \lsl_\ell(\scA)\] the
isomorphism obtained by restriction of $f$, we have an isomorphism
$\bar{f}_\scA \co \lsl_\ell(\bar{\scA}') \to \lsl_\ell(\bar \scA)$ where now
both $\bar{\scA}'$ and $\bar \scA$ are fgc quantum tori over the finite field
$R/\mmm$ of very good characteristic. This allows us to apply
Theorem~\ref{jaseco}. In view of Remark~\ref{exte-rem} we get that either
\[ (\bar{f}_\scA)_\gl \co \gl_\ell (\bar{\scA}') \to \gl_\ell(\bar \scA) \quad
\hbox{or } \quad ( \bar{f}_\scA\circ \io\op)_\gl \co
\gl_\ell(\bar{\scA}'^\op) \to \gl_\ell(\bar \scA) \] is an isomorphism of the
underlying associative algebras.

On the other side, arguing as in Lemma~\ref{exte} we get the following:
if $\euZ(\scA)$ (resp. $\euZ(\scA')$) denotes the centre
of $\scA$ (resp. $\scA'$) our map $f_{\scA}$ has a canonical extension to
\[\euZ(\scA')E_{\ell}\oplus \lsl_{\ell}(\scA')\to \euZ(\scA)E_\ell\oplus \lsl_\ell(\scA)\]
which  abusing notation we will still denote by $(f_\scA)_{\gl}$.
Note that $(f_\scA)_\gl$ coincides with the restriction of
$f_\gl$ to $\euZ(\scA')E_{\ell}\oplus \lsl_{\ell}(\scA')$
and that by our construction all matrices
$q_iE_\ell,\,q_i'E_{\ell'},\, g_i'$ and $g_i$ live in
$\euZ(\scA')E_{\ell}\oplus \lsl_{\ell}(\scA')$ and $\euZ(\scA)E_{\ell}\oplus \lsl_{\ell}(\scA)$.
Passing to the residues we get an isomorphism
\[ \overline{ (f_{\scA})_\gl} \co
\overline{\euZ(\scA')}\oplus \lsl_\ell(\bar{\scA}')
\to \overline{\euZ(\scA)}\oplus \lsl_\ell(\bar{\scA}).\]
It is easily seen from the construction that

\[\overline{(f_\scA)_\gl} =
((\bar{f}_\scA)_\gl)|_{\overline{\euZ(\scA')}\oplus\lsl(\bar{\scA}')} =:\psi.
\]
We now obtain a contradiction: In case $(\bar{f}_\scA)_{\gl}$ is
an isomorphism of the underlying associative algebras we have
\[\overline{g'_1} = \psi(\overline{g_1g_2}) = \psi(\overline{g_1}) \psi(\overline{g_2})
     = \overline{g_2'} \overline{g_3'},\]
whence \[\overline{g'_4} = \overline{g'_1 - g'_2g'_3} = \overline{g'_1}
-\overline{g_2'} \overline{g_3'} = 0\] contradicting $\overline{g'_4} \ne 0$
by Theorem~\ref{SpecializingIsomorphisms}. In the other case, one obtains a
contradiction in the same way.
\end{proof}

\subsection{Final step.} \label{step3} As indicated above, from now on we will assume  that \[ \phi:{\rm
M}_{\ell}(Q')\to {\rm M}_{\ell}(Q)\] is an isomorphism of associative
$k$-algebras. Let
$$V=Q\oplus\ldots\oplus Q$$
be the free right $Q$-module of rank $\ell$ defined in \eqref{def:lsl1} for
$\scA=Q$. We know that $\Mat_\ell(Q)$ acts on $V$ from the left while $Q$
acts from the right. We denote by $\sfB = \{e_1, \ldots, e_\ell\}$ the
standard basis of the $Q$--module $V$, defined in \eqref{def:st}.
Furthermore, we know that $E'_i=E'_{ii}\in {\rm M}_{\ell}(Q')$,
$i=1,\ldots,\ell$ form a complete orthogonal system of idempotents in
$\Mat_\ell(Q')$. Since $\phi$ preserves the associative multiplication, the
image of the standard orthogonal system $(E'_{11}, \ldots, E'_{\ell \ell})$
of $\Mat_\ell(Q')$ is a complete orthogonal system in $\Mat_\ell(Q)$. We put
\[
       \widetilde{E}_i=\phi(E'_{ii})\in \Mat_\ell(Q), 1 \le i \le \ell.\]
We then know from Lemma~\ref{comoc} that $V$ decomposes with respect to
$(\widetilde E_1, \ldots, \widetilde E_\ell)$:
\[ V = V_1 \oplus \cdots \oplus V_\ell, \quad \hbox{for}\quad V_i=\widetilde{E}_i(V).
\]
As shown in Lemma~\ref{comoc}\eqref{comoc-c} conjugacy will follow once we
know that all the $V_i$ are cyclic $Q$--modules. We will prove this using
again specialization. \sm

To simplify the notation we let $U= V_i$ for any one of the $i$, $1 \le i \le
\ell$. We can apply the results of \S\ref{sec:E} and choose a minimal vector
$u_0 \in U^+$.  For $t\in \NN$ define
\begin{align*}
  P_t &= \{ q \in Q^+: \deg(q) \le t \}, \\
  V_t &= \{ v\in V^+ : \deg(v) \le t \}, \\
  U_t &= V_t \cap U^+.
\end{align*}
The spaces $P_t$, $V_t$ and $U_t$ are finite-dimensional $k$-vector spaces.
We denote by $\bbP(\cdot)$ the corresponding projective spaces. Since $0 \ne
q \implies u_0 q \ne 0$ we have a well-defined regular map
\[
   \vphi_t \co \bbP(P_t) \to \bbP( U_{t + \deg(u_0)}), \quad [q] \mapsto [u_0 q].
\]
Its image is the zero set of a finite set $G_t$ of non-zero homogenous
polynomials (in fact linear forms) with coefficients in $k$: \[
   \mathrm{Im}(\vphi_t) = \mathrm{Zero}( G_t).
\]
Similarly, for $0 \le s < \deg(u_0)$ we have a regular map
\[
   \ga_s \co \bbP(V_s) \times \bbP( P_{\deg(u_0) - s}) \to \bbP(V_{\deg(u_0)}),
   \quad ([v], [q]) \to [vq].
\]
Since we are dealing with projective spaces, the image of $\ga_s$ is a closed
subvariety, whence given by a finite set $H_s$ of non-zero homogeneous
polynomials with coefficients in $k$:
\[
   \mathrm{Im} (\ga_s) = \mathrm{Zero}( H_s).
\]
By Lemma~\ref{lem812}, $[u_0] \not\in \mathrm{Im}(\ga_s)$ for all $0 \le s <
\deg(u_0)$. Hence:
\begin{equation} \label{step3-1}
   h_s(u_0) \ne 0 \quad \hbox{for some $h_s \in H_s$, $0 \le s <
\deg(u_0)$.}
\end{equation}
Recall that our goal is to show $U=u_0 Q$, i.e., in view of
Lemma~\ref{lem811}: $U^+ = u_0 Q^+$. For the purpose of contradiction, assume
this is not the case. Thus there exists $v_0 \in U^+ \setminus u_0 Q^+$.
Observe \[ d: = \deg(v_0) - \deg(u_0) \ge 0.\] Therefore $[v_0] \not\in
\mathrm{Im}(\vphi_d)$, i.e,
\begin{equation}   \label{step3-2}
   g_d(v_0) \ne 0 \quad \hbox{for some $g_d\in G_d$.}
\end{equation}
We now apply Corollary~\ref{fgc-spec} to construct a subring $R< k$. The
finitely many elements $a_i \in k$, $b_i \in Q$ and $g_i \in \Mat_\ell(Q)$ of
loc.\ cit.\ are the following.
\begin{enumerate}
  \item [$\bullet$] in $k$: the elements $h_s(u_0)$ and $g_d(v_0)$ of
      \eqref{step3-1} and \eqref{step3-2} respectively; all $q_{ij}$; the
      coefficients of the polynomial $g_d$ and of all polynomials in $H_s$,
      $0 \le s < \deg(u_0)$;

 \item[$\bullet$] in $Q$: the (by definition non-zero) coefficients of the
     vectors $u_0$ and $v_0$;

  \item[$\bullet$] in $\Mat_\ell(Q)$: the matrices $\widetilde E_i$, $1 \le
      i \le \ell$.
  \end{enumerate}
As in Corollary~\ref{fgc-spec} let $\scA= \bigoplus_{\la \in \La} Rx^\la$. We
then have objects ``over $R$'':
\begin{align*}
  V_\scA = \textstyle \bigoplus_{i=1}^\ell e_i \scA, \quad
   \scA^+ = \bigoplus_{\la \in \La^+} R x^\la, \quad
    V^+_\scA = \bigoplus_{i= 1}^\ell e_i \scA^+.
\end{align*}
Since all matrices $\widetilde E_i \in \Mat_\ell(\scA)$ we get a decomposition
\[
   V_\scA = V_{\scA,1} \oplus \cdots \oplus V_{\scA,\ell}, \quad V_{\scA,i} = V_\scA \cap
V_i. \] In particular $U_\scA = U \cap V_\scA$. We choose the maximal ideal
$\mmm\triangleleft R$ as in Corollary~\ref{fgc-spec}, and denote by
$\overline{\phantom{0}}$ the quotient objects:
\begin{align*}
  \bar R &= R/\mmm, \\
 \bar \scA  = \scA / \mmm \scA &= \textstyle \bigoplus_{\la \in \La} \bar R x^\la, &
 \bar \scA^+ = \scA^+ / \mmm \scA^+ &= \textstyle \bigoplus_{\la \in \La^+} \bar R x^\la,\\
 \bar V_\scA = V_\scA / \mmm V_\scA &= \textstyle \bigoplus_{i=1}^\ell e_i \bar \scA \simeq \bar \scA^\ell, &
\bar V_\scA^+ = V_\scA^+ / \mmm V_\scA^+ &= \textstyle \bigoplus_{i=1}^\ell e_i \bar \scA^+ \simeq \bar \scA^\ell,\\
 \bar U_\scA &= U_\scA / \mmm U_\scA.
\end{align*}
By construction, $\bar\scA$ is a fgc quantum torus over the finite field $\bar
R$ which has very good characteristic for $\lsl_\ell(\bar{\scA})$. Hence, by
Corollary~\ref{cor-conj}, conjugacy holds in $\lsl_\ell(\bar \scA)$. Thus, by
Lemma~\ref{comoc}, $\bar U_\scA$ is a free $\bar \scA$--module, say $\bar U_\scA =
\bar c \cdot \bar \scA$. We can apply the results of \S\ref{sec:E}: without loss
of generality, $\bar c \in \bar U_\scA^+= \bar{U}_\scA \cap \bar{V}_\scA^+$. We can even assume that
$\bar c$ is indivisible. Thus, by Corollary~\ref{genco}, $\bar U_\scA^+ = \bar c
\bar \scA^+$. Since $\bar u_0 \in \bar U_\scA^+$ we get from Lemma~\ref{lem811}
that
\[   \bar u_0 = \bar c \cdot \bar a \quad \hbox{for some $\bar a \in \bar \scA^+$}.
\]

Our next goal is to show that $\bar a \in \bar R \cdot 1_{\bar \scA}$. To this
end we use the ``bar''-versions of the vector spaces and maps defined above:
\begin{align*}
  \bar P_t &= \{ q \in \bar \scA^+: \deg(\bar q) \le t \}, \quad
  \bar V_t = \{ \bar v\in \bar V^+_\scA : \deg(\bar v) \le t \}, \quad
   \bar U_t = \bar V_t \cap \bar U^+, \\
 \bar \vphi_t &\co \bbP(\bar P_t) \to \bbP( \bar U_{t + \deg(u_0)}),
                \quad [\bar q] \mapsto [\bar u_0 \bar q], \\
  \bar \ga_s  & \co \bbP(\bar V_s) \times \bbP( \bar P_{\deg(u_0) - s})
             \to \bbP(\bar V_{\deg(u_0)}), \quad ([\bar v], [\bar q]) \to
                        [\bar v \bar q]
\end{align*}
By base change, $\mathrm{Im}( \bar \vphi_t)$ is the zero set of the
polynomials $\{\bar g: g\in G_t\}$. Similarly, $\mathrm{Im}(\bar \ga_s)$ is
the zero set of the polynomials $\bar h$, $h \in H_s$. From $\bar u_0 = \bar
c \cdot \bar a$ we obtain $\deg(\bar u_0) = \deg(\bar c) + \deg (\bar a)$.
Assuming $\deg(\bar a) >0$, it follows that $\bar u_0 \in \mathrm{Im}(\bar
\ga_s)$ for $0 \le s= \deg(\bar c) < \deg(\bar u_0) = \deg(u_0)$. Hence $\bar
h(\bar u_0) = \overline{ h(u_0)} = 0$ for all polynomials $h \in H_s$. But
this contradicts \eqref{step3-1}: $\bar h_s(u_0) \ne 0$ by construction of
$R$ and $\mmm$. Hence $\deg(\bar a) = 0$, proving that $\bar u_0$ is also a
generator of $\bar U_\scA^+$: $\bar U_\scA^+ = \bar u_0 \bar \scA^+$.

Recall the element $v_0 \in U^+ \setminus u_0 Q^+$. We have $0 \ne \bar v_0
\in \bar U_\scA^+ = \bar u_0 \bar \scA^+$. Hence $\bar g(\bar v_0) =
\overline{ g(v_0)} = 0$ for all $g\in G_d$. But this contradicts
\eqref{step3-2}: $\overline{g_d(v_0)} \ne 0$ by construction of $R$ and
$\mmm$. Thus, we have arrived at the final contradiction: There does not
exist $v_0 \in U^+ \setminus u_0 Q^+$. It follows that $U$ is indeed
generated by $u_0$.

\end{document}